%% file: InverseStefanProblemPart2_Arxiv_submission_15May_2016.TEX
\documentclass{article}
\usepackage{amsfonts,amsmath,amssymb}
\usepackage{graphicx,epsfig}
\usepackage{checkend}
\usepackage{longtable,geometry}
\newtheorem{definition}{Definition}[section]
\begin{document}
\input dibemod.0605

\input harnack.mac

\begin{center}  {\huge\textbf{On the Optimal Control of the Free Boundary Problems for the 
Second Order Parabolic Equations. II.Convergence of the Method of Finite Differences }}
\par \medskip\bigskip\end{center}
\begin{center} {\Large\textsc{Ugur G. Abdulla}}
\par \medskip\bigskip\end{center}
\begin{center} {\large\noindent \textsc{Department of Mathematics, Florida Institute of Technology, Melbourne, Florida 32901}}
\par \medskip\bigskip\end{center}
{\bf Abstract.} We develop a new variational formulation of the inverse Stefan problem, where information on the heat flux on the fixed boundary is missing and must be found along with the temperature and free boundary. We employ optimal control framework, where boundary heat flux and free boundary are components of the control vector, and optimality criteria consist of the minimization of the sum of $L_2$-norm declinations from the available measurement of the temperature flux on the fixed boundary and available 
information on the phase transition temperature on the free boundary. This approach allows 
one to tackle situations when the phase transition temperature is not known explicitly, and is available through measurement with possible error. It also allows for the development of iterative numerical methods of least computational cost due to the fact that for every given control vector, the parabolic PDE is solved in a fixed region instead of full free boundary problem. In {\it Inverse Problems and Imaging, 7, 2(2013), 307-340} we proved well-posedness in Sobolev spaces framework and 
convergence of time-discretized optimal control problems. In this paper we perform full discretization and prove convergence of the discrete optimal control problems to the original problem both with respect to cost functional and control.



\newpage
\section{Description of Main Results}\label{description of results,historical remarks}
\subsection{Introduction and Motivation}\label{E:1:1}
{\large
Consider the general one-phase Stefan problem (\cite{Friedman1, Meyrmanov}): find the temperature function $u(x,t)$ and the free boundary $x=s(t)$ from the following conditions
\begin{equation}\label{Eq:W:1:1}
(a(x,t)u_{x})_{x}+b(x,t)u_{x}+c(x,t)u-u_{t}=f(x,t),\quad \text{for}~(x,t) \in \Omega
\end{equation}
\begin{equation}\label{Eq:W:1:2}
u(x,0)=\phi(x),\qquad 0 \leq x \leq s(0)=s_{0}
\end{equation}
\begin{equation}\label{Eq:W:1:3}
a(0,t)u_{x}(0,t)=g(t),\qquad 0 \leq t \leq T
\end{equation}
\begin{equation}\label{Eq:W:1:4}
a(s(t),t)u_{x}(s(t),t) + \gamma(s(t),t)s'(t)=\chi(s(t),t),\qquad 0 \leq t \leq T
\end{equation}
\begin{equation}\label{Eq:W:1:5}
u(s(t),t)=\mu(t),\qquad 0 \leq t \leq T
\end{equation}
where $a$, $b$, $c$, $f$, $\phi$, $g$, $\gamma$, $\chi$, $\mu$ are known functions and
\begin{equation}\label{Eq:W:1:6}
a(x,t)\geq a_{0}>0,\quad s_{0}>0
\end{equation}
\[ \Omega=\left\{ (x,t) :~0<x<s(t),~0<t\leq T \right\} \]
In the physical context, $f$ characterizes the density of the sources, $\phi$ is the initial temperature, $g$ is the heat flux on the fixed boundary and $\mu$ is the phase transition temperature.  

Assume now that some of the data is not available, or involves some measurement error.  For example, assume that the heat flux $g(t)$ on the fixed boundary $x=0$ is not known and must be found along with the temperature $u(x,t)$ and the free boundary $s(t)$.  In order to do that, some additional information is needed.  Assume that this additional information is given in the form of the temperature measurement along the boundary $x=0$:
\begin{equation}\label{Eq:W:1:7}
u(0,t)=\nu(t),\quad \text{for}~0 \leq t \leq T
\end{equation}
\textbf{Inverse Stefan Problem (ISP):} \textit{Find the functions $u(x,t)$ and $s(t)$ and the boundary heat flux $g(t)$ satisfying conditions (\ref{Eq:W:1:1})-(\ref{Eq:W:1:7}}).

Motivation for this type of inverse problem arose, in particular, in the modeling of bioengineering problems on the laser ablation of biological tissues through Stefan problem (\ref{Eq:W:1:1})-(\ref{Eq:W:1:6}), where $s(t)$ is the ablation depth at the moment $t$. The boundary temperature measurement $u(0,t)$ contains an error, which makes it impossible to get reliable measurement of the boundary heat flux $g(t)$, and the ISP must be solved for its identification. This approach allows us to regularize an error contained in a measurement $\nu(t)$. Another advantage of this approach is that, in fact, condition (\ref{Eq:W:1:5}) can be treated as a measurement of the temperature on the ablation front, and our approach allows us to regularize an error contained in temperature measurement $\mu(t)$ on the ablation front. Still another important motivation arises in optimal control of the Stefan problem, where controlling $g(t)$ is equivalent of controlling external temperature along the fixed boundary. It should be pointed out that the method of this paper can be applied to different type of inverse problems. For example, (\ref{Eq:W:1:7}) can be replaced with
\[ u(x,T)=w(x), \quad \text{for}~0 \leq x \leq s(T), \]
meaning that measurements are taken for the final temperature distribution $w(x)$ and final ablation depth $s(T)$. Instead of identification of the boundary flux $g$, one can consider the inverse free boundary problem with any of the unknown coefficients $a, b, c$ or right hand side $f$.

The ISP is not well posed in the sense of Hadamard.  If there is no coordination between the input data, the exact solution may not exist.  Even if it exists, it might be not unique, and most importantly, there is no continuous dependence of the solution on the data. 
The ISP was first mentioned in \cite{Cannon3}, in the form of finding
a heat flux on the fixed boundary which provides a desired free boundary. This problem is similar to a non-characteristic Cauchy problem for the heat equation. The variational approach
for solving this ill-posed inverse Stefan problem was performed in \cite{BudakVasileva1, BudakVasileva2}. The first result on the optimal control of the Stefan problem appeared in \cite{Vasilev}. It consists of finding the optimal value of the external temperature along the fixed boundary, in order to ensure that the solutions of the Stefan problem are close to the measurements taken
at the final moment. In \cite{Vasilev}, the existence result was proved. In \cite{Yurii} the Frechet differentiability and the convergence of the difference schemes was proved for the same problem and Tikhonov regularization was suggested. Later development of the inverse Stefan problem was along these two lines: Inverse Stefan problems with given phase boundaries were considered in \cite{Alifanov,Bell,Budak,Cannon,Carasso,Ewing1,Ewing2,Hoffman,Sherman,Goldman}; optimal control of Stefan problems, or equivalently inverse problems with unknown phase boundaries were investigated in \cite{Baumeister,Fasano,Hoffman1,Hoffman2,Jochum2,Jochum1,Knabner,Lurye,Nochetto, Niezgodka,Primicero,Sagues,Talenti,Goldman}. We refer to monography \cite{Goldman} for a complete list of references of both types of inverse Stefan problems, both for linear and quasilinear parabolic equations. 
The main methods used to solve the inverse Stefan problem are based on variational formulation, method of quasi-solutions or Tikhonov regularization which takes into account ill-posedness in terms of the dependence of the solution on the inaccuracy involved in the measurement (\ref{Eq:W:1:7}), Frechet differentiability and iterative conjugate gradient methods for numerical solution. Despite its effectiveness, this approach has some deficiencies in many practical applications:
\begin{itemize}
\item Solution of the inverse Stefan problem is not continuously dependent on the phase transition temperature $\mu(t)$: small perturbation of the phase transition temperature may imply significant change of the solution to the inverse Stefan problem. Accordingly, any regularization which equally takes into account instability with respect to both $\nu(t)$ from measurement (\ref{Eq:W:1:7}), and the phase
transition temperature $\mu(t)$ from (\ref{Eq:W:1:5}) will be preferred. It should be also mentioned that in many applications the phase transition temperature is not known explicitly. In many processes the melting temperature of pure material at a given external action depends on the process evolution. For example, gallium (Ga, atomic number 31) may remain in the liquid phase at temperatures well below its mean melting temperature (\cite{Meyrmanov}). 
\item Numerical implementation of iterative gradient type methods within the existing approach requires solving the full free boundary problem at every step of 
the iteration, and accordingly has quite a high computational cost. An iterative gradient method which requires solution of the boundary value problem in a fixed region at every step  would definitely be much more effective in terms of the computational cost.  
\end{itemize}

The main goal of this project is to develop a new variational approach based on the optimal control theory which is capable of addressing both of the mentioned issues and allows the inverse Stefan problem to be solved numerically with least computational cost by using conjugate gradient methods in Hilbert spaces.
In \cite{Abdulla1} we proved the existence of the optimal control and convergence of the family of time-discretized optimal control problems to the continuous problem. In this paper we perform full discretization through finite differences and prove the convergence of the discrete optimal control problems to the continuous problem both with respect to cost functional and control.  We employ Sobolev spaces framework which allows us to reduce the regularity and structural requirements on the data. We address the problems of Frechet differentiability and application of iterative conjugate gradient methods in Hilbert spaces in an upcoming paper.

Throughout the paper we use the usual notation for Sobolev spaces according to references \cite{LSU,BIN,Nikolski,Solonnikov1,Solonnikov2}.
Notation is described below in Section~\ref{E:1:1a}.

\subsection{Notation of Sobolev Spaces}\label{E:1:1a}
$L_2[0,T]$ - Hilbert space with scalar product
\[ (u,v)=\int_0^T uvdt \]
$W_2^k[0,T], k=1,2,...$ - Hilbert space of all elements of $L_2[0,T]$ whose weak derivatives up to order $k$ belongs to $L_2[0,T]$ and scalar product is defined as
\[ (u,v)=\int_0^T \sum_{s=0}^k \frac{d^su}{dt^s} \frac{d^sv}{dt^s} dt  \]
$W_2^{\frac{1}{4}}[0,T]$ - Banach space of all elements of $L_2[0,T]$ with finite norm
\[ \Vert u\Vert_{W_2^{\frac{1}{4}}[0,T]}= \Big ( \Vert u \Vert^2_{L_2[0,T]} + \int_0^T dt \int_0^T \frac{|u(t)-u(\tau)|^2}{|t-\tau|^{\frac{3}{2}}} d\tau \Big )^{\frac{1}{2}} \]
$L_2(\Omega)$ - Hilbert space with scalar product 
\[ (u,v)=\int_{\Omega} uv dx dt  \]
$W_2^{1,0}(\Omega)$ - Hilbert space of all elements of $L_2(\Omega)$ whose weak derivative $\frac{\partial u}{\partial x}$ belongs to $L_2(\Omega)$, and scalar product is defined as
\[ (u,v)=\int_{\Omega} \Big ( uv + \frac{\partial u}{\partial x}\frac{\partial v}{\partial x} \Big ) dx dt \]
$W_2^{1,1}(\Omega)$ - Hilbert space of all elements of $L_2(\Omega)$ whose weak derivatives $\frac{\partial u}{\partial x}$, $\frac{\partial u}{\partial t}$ belong to $L_2(\Omega)$, and scalar product is defined as
\[ (u,v)=\int_{\Omega} \Big ( uv + \frac{\partial u}{\partial x}\frac{\partial v}{\partial x}+ \frac{\partial u}{\partial t}\frac{\partial v}{\partial t} \Big ) dx dt \]
$V_2(\Omega)$ - Banach space of all elements of $W_2^{1,0}(\Omega)$ with finite norm
\[ \Vert u \Vert_{V_2(\Omega)} = \Big ( esssup_{0\le t \le T} \Vert u(x,t) \Vert^2_{L_2[0,s(t)]} + \Big \Vert \frac{\partial u}{\partial x} \Big \Vert^2_{L_2(\Omega)} \Big )^{\frac{1}{2}}  \]
$V_2^{1,0}(\Omega)$ - Banach space which is the completion of $W_2^{1,1}(\Omega)$ in the norm of $V_2(\Omega)$. It consists of all elements of $V_2(\Omega)$, continuous with respect to $t$ in norm of $L_2[0,s(t)]$ and with finite norm
\[ \Vert u \Vert_{V_2^{1,0}(\Omega)} = \Big ( \max_{0\le t \le T} \Vert u(x,t) \Vert^2_{L_2[0,s(t)]} + \Big \Vert \frac{\partial u}{\partial x} \Big \Vert^2_{L_2(\Omega)} \Big )^{\frac{1}{2}}  \]
$W_2^{2,1}(\Omega)$ - Hilbert space of all elements of $L_2(\Omega)$ whose weak derivatives $\frac{\partial u}{\partial x}$, $\frac{\partial u}{\partial t}$, $\frac{\partial^2u}{\partial x^2}$ belong to $L_2(\Omega)$, and scalar product is defined as
\[ (u,v)=\int_{\Omega} \Big ( uv + \frac{\partial u}{\partial x}\frac{\partial v}{\partial x}+ \frac{\partial u}{\partial t}\frac{\partial v}{\partial t}+\frac{\partial^2 u}{\partial x^2}\frac{\partial^2 v}{\partial x^2} \Big ) dx dt \]

\subsection{Optimal Control Problem}\label{E:1:2}
Consider a minimization of the cost functional
\begin{equation}
\mathcal{J}(v)=\beta_{0}\Vert u(0,t)-\nu(t)\Vert_{L_{2}[0,T]}^{2}+\beta_{1}\Vert u(s(t),t)-\mu(t)\Vert_{L_{2}[0,T]}^2\label{Eq:W:1:8}
\end{equation}
on the control set 
\begin{gather*}
V_{R}=\{ v=(s,g) \in W_{2}^{2}[0,T]\times W_{2}^{1}[0,T]: \delta \leq s(t)\leq l, s(0)=s_0, s'(0)=0, \nonumber\\
\max (~\Vert s\Vert_{W_{2}^{2}}; ~\Vert g\Vert_{W_{2}^{1}} \leq R\}
\end{gather*}
where $\delta, l,R, \beta_0, \beta_1$ are given positive numbers, and $u=u(x,t;v)$ be a solution of the Neumann problem (\ref{Eq:W:1:1})-(\ref{Eq:W:1:4}).   

\begin{definition}
The function $u \in W_{2}^{1,1}(\Omega)$ is called a weak solution of the problem (\ref{Eq:W:1:1})-(\ref{Eq:W:1:4}) if $u(x,0)=\phi(x) \in W_{2}^{1}[0,s_0]$ and
\begin{gather}
0=\int_{0}^{T}\int_{0}^{s(t)}[ a u_{x}\Phi_{x}-bu_{x}\Phi - c u \Phi + u_{t} \Phi+f\Phi] \,dx\,dt \nonumber\\
 +\int_{0}^{T}[ \gamma(s(t),t)s'(t)-\chi(s(t),t)]\Phi(s(t),t)\, dt
+\int_{0}^{T}g(t)\Phi(0,t)\, dt\label{Eq:W:1:9}
\end{gather}
for arbitrary $\Phi \in W_{2}^{1,1}(\Omega)$
\end{definition}
We also need a notion of weak solution from $V_{2}(\Omega)$ of the Neumann problem:
\begin{definition}
The function $u \in V_{2}(\Omega)$ is called a weak solution of (\ref{Eq:W:1:1})-(\ref{Eq:W:1:4}) if
\begin{gather}
0=\int_{0}^{T}\int_{0}^{s(t)}[ a u_{x}\Phi_{x}-bu_{x}\Phi - c u \Phi - u \Phi_{t}+f\Phi] \,dx\,dt -\int_{0}^{s_0}\phi(x)\Phi(x,0)\,dx+\nonumber\\
 \int_{0}^{T}g(t)\Phi(0,t)\, dt+\int_{0}^{T}[ \gamma(s(t),t)s'(t)-u(s(t),t)s'(t)-\chi(s(t),t)]\Phi(s(t),t)\, dt \label{Eq:W:1:10}
\end{gather}
for arbitrary $\Phi \in W_{2}^{1,1}(\Omega)$ such that $\left. \Phi\right|_{t=T}=0$.
\end{definition}
If $u$ is a weak solution either from $V_2(\Omega)$ (or $W_{2}^{1,1}(\Omega)$), then traces $\left. u\right|_{x=0}$ and
$\left. u\right|_{x=s(t)}$ are elements of $L_2[0,T]$, when $s\in W_2^2[0,T]$ (\cite{Nikolski, LSU}) and cost functional $\mathcal{J}(v)$
is well defined. Furthermore, formulated optimal control problem will be called Problem $I$.
\subsection{Discrete Optimal Control Problem}\label{E:1:3}
Let
\[\omega_{\tau}=\{ t_{j}=j \cdot \tau,~j=0,1,\ldots,n\} \]
be a grid on $[0,T]$ and $\tau=\frac{T}{n}$. Consider a discretized control set
\begin{equation*}
V^n_{R}=\{ [v]_{n}=([s]_n,[g]_n) \in {\mathbb R}^{2n+2}:~0<\delta\leq s_{k} \leq l,~ \max(\Vert [s]_{n}\Vert_{w_{2}^{2}}^2; ~\Vert [g]_{n}\Vert_{w_{2}^{1}}^2) \leq R^2\}
\end{equation*}
where,
\[ [s]_n=(s_0,s_1,...,s_n) \in {\mathbb R}^{n+1}, \ [g]_n=(g_0,g_1,...,g_n) \in {\mathbb R}^{n+1} \]
\[ 
\Vert [s]_{n}\Vert_{w_{2}^{2}}^2= \sum\limits_{k=0}^{n-1}\tau s_k^2+\sum\limits_{k=1}^{n}\tau s_{\overline{t},k}^2+\sum\limits_{k=0}^{n-1}\tau s_{\overline{t}t,k}^2, \  \Vert [g]_{n}\Vert_{w_{2}^{1}}^2= \sum\limits_{k=0}^{n-1}\tau g_k^2+\sum\limits_{k=1}^{n}\tau g_{\overline{t},k}^2.
\]
where we assign $s_{-1}=s_0$ and use the standard notation for the finite differences:
\[ s_{\overline{t},k}=\frac{s_k-s_{k-1}}{\tau}, \ s_{t,k}=\frac{s_{k+1}-s_{k}}{\tau}, \ s_{\overline{t}t,k}=\frac{s_{k+1}-2s_k+s_{k-1}}{\tau^2}. \]
Introduce two mappings $\mathcal{Q}_n$ and $\mathcal{P}_n$ between continuous and discrete control sets:
\[ \mathcal{Q}_n(v)=[v]_n=([s]_n,[g]_n), \quad \text{for}~ v\in V_R \]
where $s_k=s(t_k), g_k=g(t_k), k=0,1,...,n$.
\[ \mathcal{P}_n([v]_n)=v^n=(s^n,g^n)\in W_2^2[0,T]\times W_2^1[0,T] \quad \text{for}~ [v]_n \in V_R^n, \]
where
\begin{equation}\label{Eq:W:1:11}
s^n(t)=
\left\{
\begin{array}{l}
s_0+\frac{t^2}{2\tau} s_{\overline{t},1} \ \ 0\le t \le \tau,\\
s_{k-1}+(t-t_{k-1}-\frac{\tau}{2})s_{\overline{t},k-1}+\frac{1}{2}(t-t_{k-1})^2 s_{\overline{t}t,k-1} \ \ t_{k-1}\le t \le t_k, k=\overline{2,n}.
\end{array}\right.
\end{equation}
\begin{equation*}
g^n(t)=g_{k-1}+\frac{g_k-g_{k-1}}{\tau}(t-t_{k-1}), \ \ t_{k-1} \le t \le t_k, k=\overline{1,n}.
\end{equation*}
Let us now introduce spatial grid. Given $[v]_n\in V_R^n$, let $(p_0,p_1,\cdots,p_n)$ be a permutation of $(0,1,\cdots,n)$ according to order
\[ s_{p_0}\le s_{p_1}\le \cdots \le s_{p_n} \]
In particular, according to this permutation for arbitrary $k$ there exists a unique $j_k$ such that
\begin{equation}\label{Eq:W:1:11a}
s_k=s_{p_{j_k}}
\end{equation}
Furthermore, unless it is necessary in the context, we are going to write simply $j$ instead of subscript $j_k$. Let
\[\omega_{p_0}=\{ x_{i}: x_i=i \cdot h,~i=0,1,\ldots,m_0^{(n)}\} \]
be a grid on $[0,s_{p_0}]$ and $h=\frac{s_{p_0}}{m_0^{(n)}}$. Furthermore we always assume that
\begin{equation}\label{htau} 
h = O(\sqrt{\tau}), \quad \text{as}~ \tau \rightarrow 0.
\end{equation}
We continue construction of the spatial grid by induction. Having constructed $\omega_{p_{k-1}}$ on $[0,s_{p_{k-1}}]$ we construct
\[ \omega_{p_k}=\{ x_i:~i=0,1,\cdots, m_k^{(n)} \} \]
on $[0,s_{p_{k}}]$, where $m_k^{(n)}\ge m_{k-1}^{n}$, and inequality is strict if and only if $s_{p_{k}}>s_{p_{k-1}}$; for $i\le m_{k-1}^{(n)}$ points $x_i$ are the same as in grid $\omega_{p_{k-1}}$.
Finally, if $s_{p_{n}}<l$, then we introduce a grid on $[s_{p_n},l]$
\[ \overline{\omega}=\{x_i: x_i=s_{p_n}+(i-m_n^{(n)}) \overline{h}, ~i=m_n^{(n)},\cdots, N \} \]
of stepsize order $h$, i.e. $\overline{h}=O(h)$ as $h \rightarrow 0$. 
Furthermore we simplify the notation and write $m_k^{(n)}\equiv m_k$. Let
\[ h_i=x_{i+1}-x_i, \ i=0,1,\cdots,N-1; \]
and assume that
\begin{equation*} 
m_k \rightarrow +\infty, \quad \text{as}~ n\rightarrow \infty.
\end{equation*}
Introduce Steklov averages
\[ d_{k}(x)=\frac{1}{\tau}\int_{t_{k-1}}^{t_{k}}d(x,t)\,dt, \ h_{k}=\frac{1}{\tau}\int_{t_{k-1}}^{t_{k}}h(t)\,dt, \ d_{ik}=\frac{1}{h_i \tau} \int_{x_i}^{x_{i+1}}\int_{t_{k-1}}^{t_k} d(x,t)\,dt\,dx, \]
where $i=0,1,\cdots,N-1; \ k=1,\cdots,n;$ $d$ stands for any of the functions $a$, $b$, $c$, $f$, and $h$ stands for any of the functions $\nu$, $\mu$, $g$ or $g^n$. Given $v=(s,g) \in V_R$ we define Steklov averages of traces
\begin{equation}\label{Eq:W:1:12}
\chi^{k}_s=\frac{1}{\tau} \int_{t_{k-1}}^{t_{k}}\chi(s(t),t) \,dt, \
(\gamma_s s')^k=\frac{1}{\tau} \int_{t_{k-1}}^{t_{k}}\gamma(s(t),t)s'(t) \,dt.
\end{equation}
Given $[v]_n=([s]_n,[g]_n) \in V_R^n$ we define Steklov averages $\chi^{k}_{s^n}$ and $(\gamma_{s^n} (s^n)')^k$ through (\ref{Eq:W:1:12}) with  $s$ replaced by $s^n$ from (\ref{Eq:W:1:11}).

Let $\phi^n$ be a piecewise constant approximation of $\phi$:
\[ \phi^n(x)=\phi_i:= \phi(x_i),  \quad \text{for}~ x_i< x \leq x_{i+1}, i=0,..,N-1  \]

Next we define a discrete state vector through discretization of the integral identity (\ref{Eq:W:1:9})
\begin{definition}\label{discretestatevector}
Given discrete control vector $[v]_n$, the vector function
\[ [u([v]_n)]_n=(u(0),u(1),...,u(n)), \ u(k)\in {\mathbb R}^{N+1}, \ k=0,\cdots,n  \]
is called a discrete state vector if 
\begin{description}
\item{\bf(a)} First $m_0+1$ components of the vector $u(0)\in {\mathbb R}^{N+1}$ satisfy 
\[ u_i(0)=\phi_i := \phi(x_i), \ i=0,1,\cdots,m_0; \] 
\item{\bf(b)} Recalling \eqref{Eq:W:1:11a}, for arbitrary $k=1,\cdots,n$ first $m_j+1$ components of the vector $u(k)\in {\mathbb R}^{N+1}$ solve the following system of $m_j+1$ linear algebraic equations: 
\begin{gather}
\Big [ a_{0k}+hb_{0k}-h^2c_{0k}+\frac{h^2}{\tau} \Big ] u_0(k) - \Big [ a_{0k}+hb_{0k} \Big ] u_1(k)=\frac{h^2}{\tau}u_0(k-1)-h^2f_{0k}-hg^n_{k}, \nonumber\\
-a_{i-1,k}h_iu_{i-1}(k)+\Big [ a_{i-1,k}h_i+a_{ik}h_{i-1}+b_{ik}h_ih_{i-1}-c_{ik}h_i^2h_{i-1}+\frac{h_i^2h_{i-1}}{\tau} \Big ] u_i(k)- \nonumber\\
\Big [ a_{ik}h_{i-1}+b_{ik}h_ih_{i-1} \Big ] u_{i+1}(k) = -h_i^2h_{i-1}f_{ik}+\frac{h_i^2h_{i-1}}{\tau}u_i(k-1), \ i=1,\cdots,m_j-1  \nonumber\\
-a_{m_j-1,k} u_{m_j-1}(k)+a_{m_j-1,k} u_{m_j}(k)=-h_{m_j-1} \Big [  (\gamma_{s^n} (s^n)')^k-\chi^{k}_{s^n} \Big ].\label{alma}
\end{gather}
\item{\bf(c)} For arbitrary $k=0,1,...,n$, the remaining components of $u(k)\in {\mathbb R}^{N+1}$ are calculated as
\[ u_i(k)= \hat{u}(x_i;k), \ m_j\le i \le N  \]
where $\hat{u}(x;k) \in W_2^1[0,l]$ is a piecewise linear interpolation of $\{u_i(k): \ i=0,\cdots,m_j \}$, that is to say
\[  \hat{u}(x;k)=u_i(k)+\frac{u_{i+1}(k)-u_i(k)}{h_i} (x-x_i), \ x_i\le x\le x_{i+1}, i=0,\cdots,m_j-1,  \]
iteratively continued to $[0,l]$ as
\begin{equation}\label{Eq:W:1:14}
\hat{u}(x;k)=\hat{u}(2^ns_k-x;k), \ 2^{n-1}s_k \le x \le 2^ns_k, n=\overline{1,n_k}, \ n_k\le n_*=1+\log_2\Big [ \frac{l}{\delta}\Big ]
\end{equation}
where $[r]$ means integer part of the real number $r$.
\end{description}
\end{definition}

It should be mentioned that for any $k=1,2,\cdots,n$ system \eqref{alma} is equivalent to the following summation identity
\begin{gather}
\sum_{i=0}^{m_j-1}h_i \Big [ a_{ik}u_{ix}(k)\eta_{ix}-b_{ik}u_{ix}(k)\eta_i-c_{ik}u_i(k)\eta_i+
f_{ik}\eta_i+u_{i\overline{t}}(k)\eta_i \Big ] + \nonumber\\ \Big [ (\gamma_{s^n} (s^n)')^k-\chi^{k}_{s^n} \Big ]\eta_{m_j}+g^n_k \eta_0=0,\label{Eq:W:1:19}
\end{gather}
for arbitrary numbers $\eta_i, i=0,1,\cdots,m_j$.

Consider a discrete optimal control problem of minimization of the cost functional
\begin{equation}\label{Eq:W:1:15}
\mathcal{I}_n([v]_n)=\beta_{0}\tau\sum\limits_{k=1}^n \Big ( u_0(k)-\nu_k \Big )^2+\beta_{1}\tau\sum\limits_{k=1}^n \Big ( u_{m_k}(k)-\mu_k \Big )^2
\end{equation}
on a set $V_R^n$ subject to the state vector defined in Definition 1.3. Furthermore, formulated discrete optimal control problem will be called Problem $I_n$.

Throughout, we use piecewise constant and piecewise linear interpolations of the discrete state vector:
given discrete state vector $[u([v]_n)]_n=(u(0),u(1),...,u(n))$, let
\[ u^\tau(x,t)=\hat{u}(x;k), \quad \text{if}~ t_{k-1}<t\le t_k, \ 0\le x \le l, \ k=\overline{0,n}, \]
\[ \hat{u}^\tau(x,t)=\hat{u}(x;k-1)+\hat{u}_{\overline{t}}(x;k)(t-t_{k-1}), \quad \text{if}~ t_{k-1}<t\le t_k, \ 0\le x \le l, \ k=\overline{1,n}, \]
\[ \hat{u}^\tau(x,t)= \hat{u}(x;n), \quad \text{if}~ t\ge T, \ 0\le x \le l. \]
\[ \tilde{u}^\tau(x,t)=u_i(k), \quad \text{if}~ t_{k-1}< t \le t_k, \ x_i \le x < x_{i+1}, \ k=\overline{1,n}, \ i=\overline{0,N-1}.\]
Obviously, we have
\[ u^\tau \in V_2(D), \ \ \hat{u}^\tau \in W_2^{1,1}(D), \ \ \tilde{u}^\tau \in L_2(D).  \]
As before, we employ standard notations for difference quotients of the discrete state vector:
\[ u_{ix}(k)=\frac{u_{i+1}(k)-u_i(k)}{h_i}, \ u_{i\overline{t}}=\frac{u_i(k)-u_i(k-1)}{\tau}, \ \quad \text{etc.}  \]

\subsection{Formulation of the Main Result}\label{E:1:4}
Let 
\[ D=\{ (x,t):~0<x<l,~0<t\leq T\}\]
Throughout the whole paper we assume the following conditions are satisfied by the data:
\[ a,b,c \in L_{\infty}(D), \ f \in L_2(D), \]
\[ \phi \in W_2^1[0,s_0], \ \gamma, \chi \in W_2^{1,1}(D), \ \mu,\nu \in L_2[0,T], \]
the coefficient $a$ satisfies \eqref{Eq:W:1:6} almost everywhere on $D$, 
the weak derivatives $\frac{\partial a}{\partial t},  \frac{\partial a}{\partial x}$ exists and 
\begin{equation}\label{conditionon_a}
\frac{\partial a}{\partial x}\in L_{\infty}(D), \ \ \int_{0}^{T}esssup_{0 \leq x \leq l} \left| \frac{\partial a}{\partial t}\right| dt <+\infty.
\end{equation}
Our main theorems read:
\begin{theorem}\label{existence}
The Problem $I$ has a solution, i.e.
\[ V_*=\{ v\in V_R: \ \mathcal{J}(v)=\mathcal{J}_* \equiv \inf\limits_{v\in V_R} \mathcal{J}(v) \} \neq \emptyset \]
\end{theorem}
Note that Theorem~\ref{existence} was already proved in  \cite{Abdulla1} by using method of lines.
\begin{theorem}\label{convergence}
Sequence of discrete optimal control problems $I_n$ approximates the optimal control problem $I$ with respect to functional, i.e.
\begin{equation}\label{Eq:W:1:18}
\lim\limits_{n\to +\infty} \mathcal{I}_{n_*}=\mathcal{J}_*, 
\end{equation}
where
\[ \mathcal{I}_{n_*}=\inf\limits_{V_R^n} \mathcal{I}_n([v]_n), \ n=1,2,... \]
If $[v]_{n_\ep}\in V_R^n$ is chosen such that
\[ \mathcal{I}_{n_*} \le \mathcal{I}_n([v]_{n_\ep})\le \mathcal{I}_{n_*}+\ep_n, \ \ep_n \downarrow 0, \]
then the sequence $v_n=(s_n,g_n)=\mathcal{P}_n([v]_{n_\ep})$ converges to some element $v_*=(s_*,g_*) \in V_*$ weakly
in $W_2^2[0,T] \times W_2^1[0,T]$, and strongly in $W_2^1[0,T] \times L_2[0,T]$. In particular $s_n$ converges to $s_*$ uniformly on $[0,T]$. For any $\delta>0$, define
  \[
    \Omega_*^{\prime} = \Omega_* \cap \{x<s_*(t)-\delta,~0<t<T\}
  \] Then the piecewise linear interpolation $\hat{u}^\tau$ of the 
discrete state vector $[u[v]_{n_\ep}]_n$ converges to the solution $u(x,t;v_*) \in W_2^{1,1}(\Omega_*)$ of the Neumann 
problem (\ref{Eq:W:1:1})-(\ref{Eq:W:1:4}) weakly in $W_2^{1,1}(\Omega_*^{\prime})$.
\end{theorem} 
\section{Preliminary Results}\label{preliminaries}
In Lemma~\ref{existencediscretestate} below we prove existence and uniqueness 
of the discrete state vector $[u([v]_n)]_n$ (see Definition~\ref{discretestatevector}) for arbitrary discrete control vector $[v]_n \in V_R^n$.
In Lemma~\ref{generalcriteria} we remind a general approximation criteria for the optimal control problems from   (\cite{Vasilev1}). In Lemma~\ref{mappings} we recall some properties of the mappings
$\mathcal{Q}_n$ and $\mathcal{P}_n$ between continuous and discrete control sets.

\begin{lemma}\label{existencediscretestate}
For sufficiently small time step $\tau$, there exists a unique discrete state vector $[u([v]_n)]_n$ for arbitrary discrete control vector $[v]_n \in V_R^n$.
\end{lemma}

Proof. As it is mentioned above, for any $k=1,2,\cdots,n$ system \eqref{alma} is equivalent to the summation identity \eqref{Eq:W:1:19}
for arbitrary numbers $\eta_i, i=0,1,\cdots,m_j$. Let $\{\tilde{u}_i(k)\}$ be a solution of the homogeneous system related to \eqref{alma}, i.e.
\[ g^n_k=(\gamma_{s^n} (s^n)')^k=\chi^{k}_{s^n}=f_{ik}=u_i(k-1)=0. \]
By choosing in \eqref{Eq:W:1:19} $\eta_i=\tilde{u}_i(k)$ we have
\begin{equation}\label{Eq:W:1:20}
\sum_{i=0}^{m_j-1}h_i a_{ik}\tilde{u}_{ix}^2(k)+\frac{1}{\tau}\sum_{i=0}^{m_j-1}h_i \tilde{u}_{i}^2(k)=
\sum_{i=0}^{m_j-1}h_i \Big [ b_{ik}\tilde{u}_{ix}(k)\tilde{u}_{i}(k)+c_{ik}\tilde{u}_{i}^2(k) \Big ]
\end{equation}
Using (\ref{Eq:W:1:6}) and Cauchy inequality with $\ep >0$ we derive that
\begin{equation}\label{Eq:W:2:3}
a_0\sum_{i=0}^{m_j-1}h_i \tilde{u}_{ix}^2(k)+\frac{1}{\tau}\sum_{i=0}^{m_j-1}h_i \tilde{u}_{i}^2(k) \le \frac{\epsilon M}{2}\sum_{i=0}^{m_j-1}h_i \tilde{u}_{ix}^2(k)+ \Big ( \frac{M}{2\epsilon}+M \Big ) \sum_{i=0}^{m_j-1}h_i \tilde{u}_{i}^2(k).
\end{equation}
where
\[ M=\max \Big ( ||a||_{L_{\infty}(D)}; ||b||_{L_{\infty}(D)}; ||c||_{L_{\infty}(D)} \Big ). \]
By choosing $\ep=a_0/M$ in (\ref{Eq:W:2:3}) we have
\begin{equation}\label{Eq:W:2:4}
\frac{a_0}{2}\sum_{i=0}^{m_j-1}h_i \tilde{u}_{ix}^2(k)+\Big ( \frac{1}{\tau}-\frac{1}{\tau_0}\Big ) \sum_{i=0}^{m_j-1}h_i \tilde{u}_{i}^2(k)\le 0,
\end{equation}
where
\[ \tau_0 = \Big ( \frac{M^2}{2a_0}+M \Big )^{-1}. \]
From (\ref{Eq:W:2:4}) it follows that $\tilde{u}_i(k)= 0$, $i=0,1,\cdots,m_j$, and hence the homogeneous system only has a trivial solution for $\tau < \tau_0$. Accordingly, system \label{Eq:W:1:13} is uniquely solvable and therefore, for any given discrete control vector $[v]_n$ there exists a unique discrete state vector defined by Definition 1.3. Lemma is proved.

The following known criteria will be used in the proof of Theorem~\ref{convergence}.  

\begin{lemma}\label{generalcriteria}
\cite{Vasilev1} Sequence of discrete optimal control problems $I_n$ approximates the
continuous optimal control problem $I$ if and only if the following conditions
are satisfied:
\begin{description}
\item{\bf(1)} for arbitrary sufficiently small $\ep>0$ there exists number $N_1=N_1(\ep)$ such that $\mathcal{Q}_N(v)\in V^n_R$ for all $v \in V_{R-\ep}$ and $N\ge N_1$; and for any fixed $\ep>0$ and for all $v\in V_{R-\ep}$ the following inequality is satisfied:
\begin{equation}\label{firstcondition}
\limsup\limits_{N\to \infty} \Big ( \mathcal{I}_N(\mathcal{Q}_N(v))-\mathcal{J}(v) \Big ) \le 0.
\end{equation}
\item{\bf(2)} for arbitrary sufficiently small $\ep>0$ there exists number $N_2=N_2(\ep)$ such that $\mathcal{P}_N([v]_N)\in V_{R+\ep}$ for all $[v]_N \in V^N_R$ and $N\ge N_2$; and for all $[v]_N\in V^N_R$, $N\ge 1$ the following inequality is satisfied:
\begin{equation}\label{secondcondition}
\limsup\limits_{N\to \infty} \Big ( \mathcal{J}(\mathcal{P}_N([v]_N))-\mathcal{I}_N([v]_N)  \Big ) \le 0.
\end{equation}
\item{\bf(3)} the following inequalities are satisfied:
\begin{align}
\limsup\limits_{\ep \to 0} \mathcal{J}_*(\ep) \ge \mathcal{J}_*, 
\ \ \liminf\limits_{\ep \to 0} \mathcal{J}_*(-\ep) \le \mathcal{J}_*,
\end{align}\label{thirdcondition}
where $\mathcal{J}_*(\pm\ep)=\inf\limits_{V_{R\pm \ep}}\mathcal{J}(u)$.
\end{description}
\end{lemma}
Next lemma demonstrates that the mappings $\mathcal{Q}_n$ and $\mathcal{P}_n$ introduced in Section~\ref{E:1:3} satisfy the conditions of Lemma~\ref{generalcriteria}.
\begin{lemma}\label{mappings}
(\cite{Abdulla1}) For arbitrary sufficiently small $\ep > 0$ there exists $n_{\ep}$ such that
\begin{equation}\label{Eq:W:2:15}
\mathcal{Q}_n(v) \in V_R^n, \quad \text{for all}~v\in V_{R-\ep} \quad \text{and}~ n>n_\ep.
\end{equation}
\begin{equation}\label{Eq:W:2:16}
\mathcal{P}_n([v]_n) \in V_{R+\ep}, \quad \text{for all}~[v]_n\in V_{R}^n \quad \text{and}~ n>n_\ep.
\end{equation}
\end{lemma}
Proof. Let $0<\ep<<R$, $v\in V_{R-\ep}$ and $\mathcal{Q}(v)=[v]_n=([s]_n,[g]_n)$. By applying Cauchy-Bunyakovski-Schwarz (CBS) inequality and Fubini's theorem we have
\begin{gather}
\sum\limits_{k=1}^{n-1}\tau s_{\ov{t}t,k}^2=\sum\limits_{k=1}^{n-1}\frac{1}{\tau^3}\Big [ \int\limits_{t_k}^{t_{k+1}}(s'(t)-s'(t-\tau))dt \Big ]^2 \le \frac{1}{\tau^2} \int_{\tau}^T|s'(t)-s'(t-\tau)|^2 dt \nonumber\\ \le \frac{1}{\tau} \int\limits_{\tau}^{T}dt \int\limits_{t-\tau}^t |s''(\xi)|^2d\xi
\le \int\limits_0^T|s''(t)|^2 dt, \ \sum\limits_{k=1}^{n}\tau s_{\ov{t},k}^2 \le \int\limits_0^T|s'(t)|^2 dt,\label{Eq:W:2:17}
\end{gather}
\begin{gather}
\tau s_{\ov{t}t,0}^2=\frac{1}{\tau^3}\Big [ \int\limits_{0}^{\tau}(s'(t)-s'(0))dt \Big ]^2 \leq \frac{1}{2} \int\limits_{0}^{\tau} |s''(t)|^2dt, \label{Eq:W:2:17a}
\end{gather}
\begin{gather}
\Big | \sum\limits_{k=0}^{n-1}\tau s_k^2 - \int_0^T s^2(t) dt \Big |= \Big | \sum\limits_{k=0}^{n-1}\int\limits_{t_k}^{t_{k+1}}\int\limits_t^{t_k}(s^2(\xi))'d\xi dt \Big | \le \nonumber\\
\sum\limits_{k=0}^{n-1}\int\limits_{t_k}^{t_{k+1}}\int\limits_{t_k}^t  [ s^2(\xi)+(s'(\xi))^2 \ ]d\xi dt
\le \tau \int\limits_0^T [ s^2(t)+(s'(t))^2 \ ]dt \le (R-\epsilon)^2 \tau, \label{Eq:W:2:18}
\end{gather}
\begin{gather}
\sum\limits_{k=1}^{n}\tau g_{\ov{t},k}^2 \le \int\limits_0^T|g'(t)|^2 dt, \ 
\Big | \sum\limits_{k=0}^{n-1}\tau g_k^2 - \int_0^T g^2(t) dt \Big | \le (R-\epsilon)^2 \tau,\label{Eq:W:2:19}
\end{gather}
From (\ref{Eq:W:2:17})-(\ref{Eq:W:2:19}) it follows that
\begin{equation}
\max \Big ( \Vert [s]_n \Vert^2_{w_2^2}, \Vert [g]_n \Vert^2_{w_2^1} \Big ) \le \max \Big ( \Vert s \Vert^2_{W_2^2[0,T]}, \Vert g \Vert^2_{W_2^1[0,T]} \Big )+(R-\epsilon)^2\tau+\frac{1}{2} \int\limits_{0}^{\tau} |s''(t)|^2dt.\label{Eq:W:2:20}
\end{equation}
From (\ref{Eq:W:2:20}), (\ref{Eq:W:2:15}) follows.

Let us know choose $[v]_n\in V_R^n$. We simplify the notation and assume $v=(s,g)=\mathcal{P}_n([v]_n)$. Through direct calculations we derive
\begin{equation}\label{Eq:W:2:21}
\Vert s \Vert_{W_2^2[0,T]}^2 \le \sum\limits_{k=0}^{n-1}\tau s_k^2 +\sum\limits_{k=1}^{n-1}\tau s_{\ov{t},k}^2 + \sum\limits_{k=0}^{n-1}\tau s_{\ov{t}t,k}^2 + \frac{1}{3}\tau s_{\ov{t},1}^2 + C \tau,
\end{equation}
where $C$ is independent of $\tau$. Furthermore, we use notation $C$ for all (possibly different) constants which are independent of $\tau$. By using CBS inequality we have
\begin{gather}
\frac{1}{3}\tau s_{\ov{t},1}^2 = \frac{4}{3\tau}(s(\tau)-s(0))^2 \le \frac{4}{3} \int\limits_0^\tau |s'(t)|^2dt\label{Eq:W:2:22}
\end{gather}
By applying Morrey inequality to $s'(t)$ from (\ref{Eq:W:2:22}) it follows 
\begin{gather}
\frac{1}{3}\tau s_{\ov{t},1}^2 \leq C \tau \Vert s \Vert_{W_2^2[0,T]}^2\label{Eq:W:2:22a}
\end{gather}
Since $[v]_n \in V_R^n$, from (\ref{Eq:W:2:21}),(\ref{Eq:W:2:22a}) it follows that for all $\tau \le (2C)^{-1}$
\begin{equation}\label{Eq:W:2:23}
\Vert s \Vert_{W_2^2[0,T]}^2 \le C,
\end{equation}
Therefore from (\ref{Eq:W:2:21}),(\ref{Eq:W:2:22a}),(\ref{Eq:W:2:23}) it follows that for sufficiently small $\tau$
\begin{equation}\label{Eq:W:2:24}
\Vert s \Vert_{W_2^2[0,T]}^2 \le \sum\limits_{k=0}^{n-1}\tau s_k^2 +\sum\limits_{k=1}^{n-1}\tau s_{\ov{t},k}^2 + \sum\limits_{k=0}^{n-1}\tau s_{\ov{t}t,k}^2 +C\tau
\end{equation}
In a similar way we calculate
\begin{equation}\label{Eq:W:2:25}
\Vert g \Vert_{W_2^1[0,T]}^2 \le \sum\limits_{k=0}^{n-1}\tau g_k^2 + \sum\limits_{k=1}^{n}\tau g_{\ov{t},k}^2 + C\tau.
\end{equation}
Hence, from (\ref{Eq:W:2:24}),(\ref{Eq:W:2:25}) it follows that for sufficiently small $\tau$
\begin{gather}
\max \Big ( \Vert s \Vert^2_{W_2^2[0,T]}, \Vert g \Vert^2_{W_2^1[0,T]} \Big ) \le \max \Big ( \Vert [s]_n \Vert^2_{w_2^2}, \Vert [g]_n \Vert^2_{w_2^1} \Big )+C\tau 
,\label{Eq:W:2:26}
\end{gather}
From (\ref{Eq:W:2:26}), (\ref{Eq:W:2:16}) follows. Lemma is proved.
\begin{corollary}\label{lipshitz}
(\cite{Abdulla1}) Let either $[v]_n \in V_R^n$ or $[v]_n = {\mathcal Q}_n(v)$ for $v \in V_R$. Then
\begin{equation}\label{Eq:W:2:28}
|s_k-s_{k-1}|\le C^\prime \tau, \ \ k=1,2,\cdots,n
\end{equation}
where $C^\prime$ is independent of $n$. 
\end{corollary}
Indeed, if $v \in V_R$, then $s^{\prime}\in W_2^1[0,T]$ and by Morrey inequality
\begin{equation}\label{Eq:W:2:29}
\Vert s^{\prime} \Vert_{C[0,T]} \le C_1 \Vert s^{\prime} \Vert_{W_2^1[0,T]} \le C_1 R
\end{equation}
and hence for the first component $[s]_n$ of $[v]_n = {\cal Q}_n(v)$ we have (\ref{Eq:W:2:28}). Also, if $[v]_n \in V_R^n$, then the sequence $v^n={\cal P}_n([v]_n)$ belongs to $V_{R+1}$ by Lemma~\ref{mappings} and 
the component $s^n$ of $v^n$ satisfies (\ref{Eq:W:2:29}). Since, $(s^n)'(t_k)=s_{\bar{t},k}, k=1,...,n$,
from (\ref{Eq:W:2:29}), (\ref{Eq:W:2:28}) follows.

Note that for the step size $h_i$ we have one of the three possibilities: $h_i=h$, or $h_i=\overline{h}$, or $h_i \le |s_{k}-s_{k-1}|$ for some $k$. Hence, from (\ref{htau}) and (\ref{Eq:W:2:28}), it follows that 
\begin{equation}\label{Eq:W:2:29a}
\max_{0\le i \le N-1}h_i = O(\sqrt{\tau}), \quad \text{as}~ \tau \rightarrow 0.
\end{equation}
\section{Proofs of the Main Results}\label{proofs}
\subsection{First Energy Estimate and its Consequences}\label{firstenergyestimate}
The main goal of this section to prove the following energy estimation for the discrete state vector.
\begin{theorem}\label{firstenergy}
For all sufficiently small $\tau$ discrete state vector $[u([v]_n)]_n$ satisfies the following stability estimations:
\begin{gather}
\max\limits_{0\le k \le n} \sum_{i=0}^{N-1} h_iu_i^2(k)+\sum_{k=1}^{n}\tau\sum_{i=0}^{N-1}h_iu_{ix}^2(k) \le \nonumber\\
C \Big ( \Vert \phi^n \Vert_{L_2(0,s_0)}^2 + \Vert g^n \Vert_{L_2(0,T)}^2 + \Vert f \Vert_{L_2(D)}^2 + 
\Vert \gamma(s^n(t),t)(s^n)'(t) \Vert_{L_2(0,T)}^2 \nonumber\\ + \Vert \chi(s^n(t),t) \Vert_{L_2(0,T)}^2 +  \sum_{k=1}^{n-1}{\bf 1}_{+}(s_{k+1}-s_{k}) \sum_{i=m_j}^{m_{j_{k+1}}-1} h_iu_i^2(k) \Big ),\label{firstenergyestimate1} 
\end{gather}
where $C$ is independent of $\tau$ and ${\bf 1}_{+}$ be an indicator function of the positive semiaxis.
\end{theorem}
First we prove the following lemma.
\begin{lemma}\label{firstenergy-step1}
For all sufficiently small $\tau$, discrete state vector $[u([v]_n)]_n$ satisfies the following estimation:
\begin{gather}
\max\limits_{1\le k \le n} \sum_{i=0}^{m_j-1} h_iu_i^2(k)+\sum_{k=1}^{n}\tau\sum_{i=0}^{m_j-1}h_iu_{ix}^2(k) + \sum_{k=1}^{n} \tau^2 \sum_{i=0}^{m_j-1}h_iu_{i\overline{t}}^2(k) \le \nonumber\\
C \Big ( \Vert \phi^n \Vert_{L_2(0,s_0)}^2 + \Vert g^n \Vert_{L_2(0,T)}^2 + \Vert f \Vert_{L_2(D)}^2 + 
\Vert \gamma(s^n(t),t)(s^n)'(t) \Vert_{L_2(0,T)}^2 \nonumber\\ + \Vert \chi(s^n(t),t) \Vert_{L_2(0,T)}^2 +  \sum_{k=1}^{n-1}{\bf 1}_{+}(s_{k+1}-s_{k}) \sum_{i=m_j}^{m_{j_{k+1}}-1} h_iu_i^2(k) \Big ),\label{firstenergyestimate3} 
\end{gather}
where $C$ is independent of $\tau$.
\end{lemma}

Proof. By choosing $\eta_i=2\tau u_i(k)$ in (\ref{Eq:W:1:19}) and by using the equality
\[ 2\tau u_{i\overline{t}}(k)u_i(k) = u_i^2(k)-u_i^2(k-1)+\tau^2 u_{i\overline{t}}^2(k)  \]
we have
\begin{gather}
\sum_{i=0}^{m_j-1}h_iu_i^2(k)-\sum_{i=0}^{m_j-1}h_iu_i^2(k-1)+\tau^2\sum_{i=0}^{m_j-1}h_iu_{i\overline{t}}^2(k)+ 2\tau\sum_{i=0}^{m_j-1}h_ia_{ik}u_{ix}^2(k) = \nonumber\\
2\tau\sum_{i=0}^{m_j-1}h_i \Big [  b_{ik}u_{ix}(k)u_i(k)+ c_{ik}u_i^2(k)
 - f_{ik}u_i(k) \Big] - \nonumber\\
 2 \tau \left[ (\gamma_{s^n} (s^n)')^k-\chi^{k}_{s^n} \right] u_{m_j}(k)
 -2\tau g^n_{k}u_0(k).\label{Eq:W:3:4}
\end{gather}
Using (\ref{Eq:W:1:6}), Cauchy inequalities with appropriately chosen  $\ep >0$, and Morrey inequality 
\begin{equation}\label{Morrey}
\max_{0\le i\le m_j}u_i^2(k) \le C_*\Vert \hat{u}(x;k) \Vert^2_{W_2^1[0,s_k]} \le C \sum_{i=0}^{m_j-1}h_i(u_i^2(k)+u_{ix}^2(k))
\end{equation}
where $C_*,C$ are independent of $\tau$ and $[u([v]_n)]_n$, from (\ref{Eq:W:3:4}) we derive that
\begin{gather}
\sum_{i=0}^{m_j-1}h_iu_i^2(k)-\sum_{i=0}^{m_j-1}h_iu_i^2(k-1)+a_0\tau\sum_{i=0}^{m_j-1}h_iu_{ix}^2(k)+\tau^2\sum_{i=0}^{m_j-1}h_iu_{i\overline{t}}^2(k)\le\nonumber\\
C_1 \tau \Big[ |(\gamma_{s^n} (s^n)')^k|^2+|\chi^{k}_{s^n}|^2+|g^n_k|^2+\sum_{i=0}^{m_j-1}h_if_{ik}^2+\sum_{i=0}^{m_j-1}h_iu_i^2(k)\Big ].\label{Eq:W:3:5}
\end{gather}
where $C_1$ is independent of $\tau$. Assuming that $\tau < C_1$, from (\ref{Eq:W:3:5}) it follows that
\begin{gather}
(1-C_1\tau)\sum_{i=0}^{m_j-1}h_iu_i^2(k) \le 
\sum_{i=0}^{m_{j_{k-1}-1}}h_iu_i^2(k-1)+ {\bf 1}_{+}(s_k-s_{k-1})\sum_{i=m_{j_{k-1}}}^{m_j-1}h_iu_i^2(k-1) + \nonumber\\ 
 C_1 \tau \left[ |(\gamma_{s^n} (s^n)')^k|^2+|\chi^{k}_{s^n}|^2+|g^n_k|^2+ \sum_{i=0}^{m_j-1}h_if_{ik}^2\right],\label{Eq:W:3:6}
\end{gather}
By induction we have
\begin{gather}
\sum_{i=0}^{m_j-1}h_iu_i^2(k) \le (1-C_1\tau)^{-k}\sum_{i=0}^{m_{j_0}-1}h_iu_i^2(0) +  \sum_{l=1}^{k} (1-C_1\tau)^{-k+l-1} \Big \{ C_1 \tau 
 \Big [ |(\gamma_{s^n} (s^n)')^l|^2+\nonumber\\ |\chi^{l}_{s^n}|^2+|g^n_l|^2+ \sum_{i=0}^{m_{j_l}-1}h_if_{il}^2 \Big ]+{\bf 1}_{+}(s_l-s_{l-1})\sum_{i=m_{j_{l-1}}}^{m_{j_l}-1}h_iu_i^2(l-1) \Big \}.\label{Eq:W:3:7}
\end{gather}
For arbitrary $1\le l \le k \le n$ we have
\begin{equation}\label{Eq:W:3.8}
(1-C_1\tau)^{-k+l-1}\le(1-C_1\tau)^{-k}\le(1-C_1\tau)^{-n}=\Big ( 1-\frac{C_1 T}{n}\Big )^{-n} \to e^{C_1T},
\end{equation}
as $\tau \to 0$. Accordingly for sufficiently small $\tau$ we have
\begin{equation}\label{Eq:W:3:9}
(1-C_1\tau)^{-k+l-1}\le 2e^{C_1T} \  \quad \text{for} \ 1\le l \le k \le n,
\end{equation}
By applying Cauchy-Bunyakovski-Schwartz (CBS) inequality from (\ref{Eq:W:3:7})-(\ref{Eq:W:3:9}) it follows that
\begin{gather}
\max\limits_{1\le k \le n} \sum_{i=0}^{m_j-1}h_iu_i^2(k) \le 
C_2 \Big ( \Vert \phi^n \Vert_{L_2(0,s_0)}^2 + \Vert g^n \Vert_{L_2(0,T)}^2  + 
\Vert \gamma(s^n(t),t)(s^n)'(t) \Vert_{L_2(0,T)}^2 + \nonumber\\ \Vert \chi(s^n(t),t) \Vert_{L_2(0,T)}^2 + \Vert f \Vert_{L_2(D)}^2 +  \sum_{l=1}^{n-1} {\bf 1}_{+}(s_{l+1}-s_{l})\sum_{i=m_{j_{l}}}^{m_{j_{l+1}}-1}h_iu_i^2(l)  \Big ).\label{Eq:W:3:10} 
\end{gather}
where $C_2$ is independent of $\tau$. Having (\ref{Eq:W:3:10}), we perform summation of (\ref{Eq:W:3:5})
with respect to $k$ from $1$ to $n$ and derive 
\begin{gather}
\sum_{i=0}^{m_{j_n}-1}h_iu_i^2(n)+a_0\sum_{k=1}^{n}\tau\sum_{i=0}^{m_j-1}h_iu_{ix}^2(k) + \sum_{k=1}^{n} \tau^2 \sum_{i=0}^{m_j-1}h_iu_{i\overline{t}}^2(k) \le \nonumber\\
\Vert \phi^n \Vert_{L_2(0,s_0)}^2 + C_3 \Big ( \Vert g^n \Vert_{L_2(0,T)}^2 + \Vert f \Vert_{L_2(D)}^2 + 
\Vert \gamma(s^n(t),t)(s^n)'(t) \Vert_{L_2(0,T)}^2 \nonumber\\ + \Vert \chi(s^n(t),t) \Vert_{L_2(0,T)}^2 +  \sum_{k=1}^{n}\tau \sum_{i=0}^{m_j-1}h_iu_i^2(k) \Big ) +\sum_{k=1}^{n-1}{\bf 1}_{+}(s_{k+1}-s_k)\sum_{i=m_{j_{k}}}^{m_{j_{k+1}}-1}h_iu_i^2(k) \label{Eq:W:3:11} 
\end{gather}
From (\ref{Eq:W:3:10}) and (\ref{Eq:W:3:11}), (\ref{firstenergyestimate3}) follows. Lemma is proved.

Proof of Theorem~\ref{firstenergy}: Due to \eqref{firstenergyestimate3}, it is enough to show that the
left hand side of \eqref{firstenergyestimate1} is bounded by the left hand side of \eqref{firstenergyestimate3}. By using reflective continuation property of $\hat{u}(x;k)$ we easily derive that
\begin{gather}
\sum_{k=1}^{n}\tau\sum_{i=0}^{N-1}h_iu_{ix}^2(k)=\tau \sum_{k=1}^{n}\int_0^l \Big | \frac{d\hat{u}(x;k)}{dx} \Big |^2 dx \le \nonumber\\ 2^{n_*}\tau \sum_{k=1}^{n}\int_0^{s_k} \Big | \frac{du(x;k)}{dx} \Big |^2 dx =2^{n_*}\sum_{k=1}^{n}\tau\sum_{i=0}^{m_j-1}h_iu_{ix}^2(k).\label{Eq:W:3:12}
\end{gather}
By using \eqref{htau} and \eqref{Eq:W:2:29} we have
\begin{gather}
\sum_{i=0}^{N-1}h_iu_i^2(k)\le 2\int_{0}^{l} \hat{u}^2(x;k) \,dx + \frac{2}{3}\sum_{i=0}^{N-1}h_i^3u_{ix}^2(k)
\le 2^{n_*+1}\int_{0}^{s_k} \hat{u}^2(x;k) \,dx +\nonumber\\ C_1 \tau \sum_{i=0}^{N-1}h_iu_{ix}^2(k) \le
2^{n_*+2}\sum_{i=0}^{m_j-1}h_iu_i^2(k)+2^{n_*+2}\sum_{i=0}^{m_j-1}\frac{1}{3}h_i^3u_{ix}^2(k) + \nonumber\\
C_1 \tau \sum_{i=0}^{N-1}h_iu_{ix}^2(k) \le 2^{n_*+2}\sum_{i=0}^{m_j-1}h_iu_i^2(k) + C_2 \tau \sum_{i=0}^{N-1}h_iu_{ix}^2(k). \label{Eq:W:3:13}
\end{gather}
From (\ref{Eq:W:3:12}),(\ref{Eq:W:3:13}) and \eqref{firstenergyestimate3}, \eqref{firstenergyestimate1} follows. Theorem is proved.

Let $[v]_n \in V_R^n, n=1,2,...$ be a sequence of discrete controls. From  Lemma~\ref{mappings}
it follows that the sequence $\{\mathcal{P}_n([v]_n)\}$ is weakly precompact in $W_2^2[0,T]\times W_2^1[0,T]$. Assume that the whole sequence converges to $v=(s,g)$
weakly in $W_2^2[0,T]\times W_2^1[0,T]$. This implies strong convergence in $W_2^1[0,T] \times L_2[0,T]$. Conversely, given control $v=(s,g)\in V_R^n$ we can choose a sequence of discrete
controls $[v]_n = \mathcal{Q}_n(v)$. Appplying Lemma~\ref{mappings} twice one can easily establish
that the sequence $\{\mathcal{P}_n([v]_n\}$ converges to $v=(s,g)$
weakly in $W_2^2[0,T]\times W_2^1[0,T]$, and strongly in $W_2^1[0,T]\times L_2[0,T]$. In the next theorem we prove the continuous dependence of the family of interpolarions $\{u^\tau\}$ on this convergence.
\begin{theorem}\label{continuity1}
Let $[v]_n \in V_R^n, n=1,2,...$ be a sequence of discrete controls and the sequence $\{\mathcal{P}_n([v]_n\}$ converges strongly in $W_2^1[0,T] \times L_2[0,T]$ to $v=(s,g)$. Then the sequence $\{u^\tau\}$ converges as $\tau \to 0$ weakly in $W_2^{1,0}(\Omega)$ to weak solution $u \in V_2^{1,0}(\Omega)$ of the problem (\ref{Eq:W:1:1})-(\ref{Eq:W:1:4}), i.e. to the solution of the integral identity (\ref{Eq:W:1:10}). Moreover, $u$ satisfies the energy estimate
\begin{equation}\label{V210estimate}
\Vert u \Vert_{V_2^{1,0}(D)}^2 \le C \Big ( \Vert \phi \Vert_{L_2(0,s_0)}^2 + \Vert g \Vert_{L_2(0,T)}^2 + \Vert f \Vert_{L_2(D)}^2 + 
\Vert \gamma \Vert_{W_2^{1,0}(D)}^2 + \Vert \chi \Vert_{W_2^{1,0}(D)}^2 \Big )
\end{equation}
\end{theorem}

Proof. In addition to quadratic interpolation of $[s]_n$ from (\ref{Eq:W:1:11}), consider two linear interpolations:
\begin{equation*}
\tilde{s}^n(t)=s_{k-1}+\frac{s_k-s_{k-1}}{\tau}(t-t_{k-1}), \ \ t_{k-1} \le t \le t_k, k=\overline{1,n}; \ \tilde{s}^n(t) \equiv s_n, \ \ t\ge T;
\end{equation*}
\begin{equation*}
\tilde{s}_1^n(t)=\tilde{s}^n(t+\tau), \ \ 0\le t \le T.
\end{equation*}
It can be easily proved that both sequences $\tilde{s}^n$ and $\tilde{s}_1^n$ are equivalent 
to the sequence $s^n$ in $W_2^1[0,T]$ and converge to $s$ strongly in $W_2^1[0,T]$. In particular,
\begin{equation}\label{Eq:W:3:19}
\sup\limits_{n}\Vert \tilde{s}_1^n \Vert_{W_2^1[0,T]} < C_*
\end{equation}
where $C_*$ is independent of $n$. 

Our next goal is to absorb the last term on the right hand side of (\ref{firstenergyestimate1}) into the left hand side. We have
\begin{gather}
\sum_{k=1}^{n-1}{\bf 1}_{+}(s_{k+1}-s_{k}) \sum_{i=m_j}^{m_{j_{k+1}}-1} h_iu_i^2(k) \le \nonumber\\ 
2\sum_{k=1}^{n-1}{\bf 1}_{+}(s_{k+1}-s_{k})\int_{s_{k}}^{s_{k+1}} \hat{u}^2(x;k) dx + \frac{2}{3} \sum_{k=1}^{n-1}{\bf 1}_{+}(s_{k+1}-s_{k}) 
 \sum_{i=m_j}^{m_{j_{k+1}}-1} h_i^3u_{ix}^2(k)\label{Eq:W:3:20}
\end{gather}
Note that if $s_{k+1}>s_k$, then all the factors $h_i$ in the second term are bounded by $s_{k+1}-s_k$ and by using  \eqref{Eq:W:2:28} we have
\begin{gather}
 \sum_{k=1}^{n-1}{\bf 1}_{+}(s_{k+1}-s_{k}) \sum_{i=m_j}^{m_{j_{k+1}}-1} h_i^3u_{ix}^2(k)\le 
(C^\prime)^2\tau \sum_{k=1}^{n-1} \tau \int_{s_k}^{s_{k+1}} \Big | \frac{d\hat{u}}{dx} \Big |^2 dx
\label{Eq:W:3:21}
\end{gather}
Due to reflective continuation we have
\begin{equation}\label{Eq:W:3:22}
\int_{s_k}^{s_{k+1}} \Big | \frac{d\hat{u}}{dx} \Big |^2 dx \le 2^{n_*-1}\int_{0}^{s_{k}} \Big | \frac{d\hat{u}}{dx} \Big |^2 dx = 2^{n_*-1}\sum_{i=0}^{m_{j}-1} h_iu_{ix}^2(k).
\end{equation}
From \eqref{Eq:W:3:21} and \eqref{Eq:W:3:22} it follows that
\begin{gather}
 \sum_{k=1}^{n-1}{\bf 1}_{+}(s_{k+1}-s_{k}) \sum_{i=m_j}^{m_{j_{k+1}}-1} h_i^3u_{ix}^2(k)\le 
2^{n_*-1}(C^\prime)^2\tau \sum_{k=1}^{n-1} \tau \sum_{i=0}^{m_{j}-1} h_iu_{ix}^2(k)
\label{Eq:W:3:23}
\end{gather}
Assuming that $\tau$ is sufficiently small and by using \eqref{Eq:W:3:20} - \eqref{Eq:W:3:23} in \eqref{firstenergyestimate1}, we absorb the last term on the right hand side of \eqref{Eq:W:3:23} into the left hand side of \eqref{firstenergyestimate1} and derive modified \eqref{firstenergyestimate1} with a new constant $C$:
\begin{gather}
\max\limits_{0\le k \le n} \sum_{i=0}^{N-1} h_iu_i^2(k)+\sum_{k=1}^{n}\tau\sum_{i=0}^{N-1}h_iu_{ix}^2(k) \le \nonumber\\
C \Big ( \Vert \phi^n \Vert_{L_2(0,s_0)}^2 + \Vert g^n \Vert_{L_2(0,T)}^2 + \Vert f \Vert_{L_2(D)}^2 + 
\Vert \gamma(s^n(t),t)(s^n)'(t) \Vert_{L_2(0,T)}^2 \nonumber\\ + \Vert \chi(s^n(t),t) \Vert_{L_2(0,T)}^2 + \sum_{k=1}^{n-1}{\bf 1}_{+}(s_{k+1}-s_{k})\int_{s_{k}}^{s_{k+1}} \hat{u}^2(x;k) dx  \Big ),\label{Eq:W:3:23e} 
\end{gather}
We can now estimate the last term on the right hand side of \eqref{Eq:W:3:23e} as in \cite{Abdulla1}:
\begin{gather}
\sum_{k=1}^{n-1}{\bf 1}_{+}(s_{k+1}-s_{k})  \int_{s_{k}}^{s_{k+1}} \hat{u}^2(x;k) dx =
\sum_{k=1}^{n-1}{\bf 1}_{+}(s_{k+1}-s_{k})  \int_{t_{k}}^{t_{k+1}}(\tilde{s}^n)^{'}(t) \hat{u}^2(\tilde{s}^n(t);k) dt = \nonumber\\ \sum_{k=1}^{n-1}{\bf 1}_{+}(s_{k+1}-s_{k})  \int_{t_{k}}^{t_{k+1}}(\tilde{s}^n)^{'}(t) \Big ( u^{\tau}(\tilde{s}^n(t),t-\tau) \Big )^2 dt = \nonumber\\ \sum_{k=1}^{n-1}{\bf 1}_{+}(s_{k+1}-s_{k})  \int_{t_{k-1}}^{t_{k}}(\tilde{s}_1^n)^{'}(t) \Big ( u^{\tau}(\tilde{s}_1^n(t),t) \Big )^2 dt. \label{Eq:W:3:24}
\end{gather}
By applying CBS inequality we have
\begin{equation}\label{Eq:W:3:25}
\Big | \sum_{k=1}^{n-1}{\bf 1}_{+}(s_{k+1}-s_{k})  \int_{s_{k}}^{s_{k+1}} \hat{u}^2(x;k) dx \Big | \le
\Vert (\tilde{s}_1^n)^{\prime}\Vert_{L_2[0,T]} \Vert u^{\tau}(\tilde{s}_1^n(t),t)\Vert_{L_4[0,T]}^2.
\end{equation}
From the results on traces of the elements of space $V_2(D)$ (\cite{LSU,BIN,Nikolski}) it follows that for arbitrary $u\in V_2(D)$ the following inequality is valid
\begin{equation}\label{Eq:W:3:26}
\Vert u(\tilde{s}_1^n(t),t)\Vert_{L_4[0,T]}\le \tilde{C}\Vert u \Vert_{V_2(D)},
\end{equation}
with the constant $\tilde{C}$ being independent of $u$ as well as $n$. From (\ref{Eq:W:3:19}),(\ref{Eq:W:3:25}) and (\ref{Eq:W:3:26}) it follows that
\begin{equation}\label{Eq:W:3:27}
\Big | \sum_{k=1}^{n-1}{\bf 1}_{+}(s_{k+1}-s_{k})  \int_{s_{k}}^{s_{k+1}} u^2(x;k) dx \Big | \le
C_* \tilde{C} \Vert u^{\tau} \Vert_{V_2(D)}^2.
\end{equation}
If the constant $C_*$ from (\ref{Eq:W:3:19}) satisfies the condition
\begin{equation}\label{Eq:W:3:28}
C_*<(C\tilde{C})^{-1}
\end{equation}
then from (\ref{Eq:W:3:23e}) and (\ref{Eq:W:3:27}) it follows that
\begin{gather}
\Vert u^{\tau} \Vert_{V_2^{1,0}(D)}^2 \le C \Big ( \Vert \phi^n \Vert_{L_2(0,s_0)}^2 + \Vert g^n \Vert_{L_2(0,T)}^2 + \Vert f \Vert_{L_2(D)}^2 + \nonumber\\
\Vert \gamma(s^n(t),t)(s^n)'(t) \Vert_{L_2(0,T)}^2 + \Vert \chi(s^n(t),t) \Vert_{L_2(0,T)}^2 \Big ),\label{Eq:W:3:29} 
\end{gather}
where $C$ is another constant independent of $n$. By applying the results on the traces of elements of $W_2^{1,0}(D)$ (\cite{BIN,Nikolski}) on smooth curve $x=s^n(t)$, Morrey inequality for $(s^n)^{\prime}$ and (\ref{Eq:W:2:16}) we have
\begin{gather}
\Vert \gamma(s^n(t),t)(s^n)'(t) \Vert_{L_2(0,T)} \le \Vert (s^n)^{\prime} \Vert_{C[0,T]} \Vert \gamma(s^n(t),t) \Vert_{L_2[0,T]} \le C_3 \Vert \gamma \Vert_{W_2^{1,0}(D)}\nonumber\\
\Vert \chi(s^n(t),t)\Vert_{L_2[0,T]} \le C_3 \Vert \chi \Vert_{W_2^{1,0}(D)}, \label{Eq:W:3:30}
\end{gather}
where $C_3$ is independent of $\gamma, \chi$ and $n$. Hence, from (\ref{Eq:W:3:29}), (\ref{Eq:W:3:30})
it follows the estimation
\begin{equation}\label{Eq:W:3:31}
\Vert u^{\tau} \Vert_{V_2^{1,0}(D)}^2 \le C \Big ( \Vert \phi^n \Vert_{L_2(0,s_0)}^2 + \Vert g^n \Vert_{L_2(0,T)}^2 + \Vert f \Vert_{L_2(D)}^2 + 
\Vert \gamma \Vert_{W_2^{1,0}(D)}^2 + \Vert \chi \Vert_{W_2^{1,0}(D)}^2 \Big ),
\end{equation}
with $C$ being independent of $n$.


If (\ref{Eq:W:3:28}) is not satisfied, then due to \eqref{Eq:W:2:28} we can partition $[0,T]$ into finitely many segments $[t_{n_{j-1}},t_{n_{j}}]$, $j=\overline{1,q}$ with $t_{n_{0}}=0$, $t_{n_{q}}=T$ in such a way that by replacing $[0,T]$ with any of the subsegments $[t_{n_{j-1}},t_{n_{j}}]$ (\ref{Eq:W:3:19}) will be satisfied with $C_*$ small enough to obey (\ref{Eq:W:3:28}). Hence, we divide $D$ into finitely many subsets
\[ D^{j}=D \cap \{ t_{n_{j-1}}<t\leq t_{n_{j}}\}\]
such that every norm $\Vert u^{\tau}\Vert_{V_{2}(D^{j})}^2$ is uniformly bounded through the right-hand side of (\ref{Eq:W:3:31}). Summation with $j=1,\ldots,q$ implies (\ref{Eq:W:3:31}).

Since $(\ph^n,g^n)$ converge to $(\phi,g)$ strongly in $L_2[0,s_0] \times L_2[0,T]$, from (\ref{Eq:W:3:31}) it follows that the sequence $\{ u^{\tau}\}$ is weakly precompact in $W_{2}^{1,0}(D)$.  Let $u \in W_{2}^{1,0}(D)$ be a weak limit point of $u^{\tau}$ in $W_{2}^{1,0}(D)$, and assume that whole sequence $\{ u^{\tau}\}$ converges to $u$ weakly in $W_{2}^{1,0}(D)$.  Let us prove that in fact $u$ satisfies the integral identity (\ref{Eq:W:1:10}) for arbitrary test function $\Phi \in W_{2}^{1,1}(\Omega)$ such that $\left. \Phi\right|_{t=T}=0$. Due to density of $C^1(\overline{\Omega})$ in $W_{2}^{1,1}(\Omega)$ it is enough to assume $\Phi \in C^1(\overline{\Omega})$. Without loss of generality we can also assume that
$\Phi \in C^1(\overline{D}_{T+\tau}), \ \Phi\equiv 0, \quad\text{for}~ T\le t \le T+\tau$, where
\[ D_{T+\tau}=\{ (x,t):~0<x<l+1,~0<t\leq T+\tau\} \]
Otherwise, we can continue $\Phi$ to $D_{T+\tau}$ with the described properties. Let 
\[ \Phi_i(k)=\Phi(x_i,t_k), \ k=0,\cdots,n+1, \ i=0,\cdots,N \]
and 
\[ \Phi^\tau(x,t)=\Phi_i(k), \Phi_x^\tau(x,t)=\Phi_{ix}(k), \Phi_t^\tau(x,t)=\Phi_{i\bar{t}}(k+1), \quad \text{for} \ t_{k-1}<t\le t_k, x_i \le x < x_{i+1}. \]
Obviously, the sequences $\{ \Phi^{\tau}\}$, $\{ \Phi_x^\tau \}$ and $\{ \Phi_t^\tau \}$ converge as $\tau \to 0$ uniformly in $\overline{D}$ to $\Phi$, $\frac{\partial \Phi}{\partial x}$ and $\frac{\partial \Phi}{\partial t}$ respectively. 
By choosing in (\ref{Eq:W:1:19}) $\eta_i=\tau \Phi_i(k)$, after summation with respect to $k=\overline{1,n}$ and transformation of the time difference term as follows
\begin{gather}
\sum_{k=1}^{n}\tau\sum_{i=0}^{m_j-1}h_iu_{i\bar{t}}(k)\Phi_i(k)=
-\sum_{k=1}^{n-1}\tau\sum_{i=0}^{m_{j_{k+1}}-1}h_iu_i(k)\Phi_{i\bar{t}}(k+1)- \sum_{i=0}^{m_{j_1}-1}h_iu_i(0)\Phi_i(1) +  \nonumber\\ \sum_{k=1}^{n-1}sign(s_k-s_{k+1})\sum_{i=\alpha_k}^{\beta_k-1}h_iu_i(k)\Phi_i(k)=
-\sum_{k=1}^{n}\int_{t_{k-1}}^{t_{k}}\int_{0}^{s_{k+1}}\tilde{u}^{\tau}\Phi_{t}^{\tau}\, dx\, dt
-\nonumber\\ \int_{0}^{s_{1}}\phi^n(x)\Phi^{\tau}(x,\tau)\, dx 
+\sum_{k=1}^{n-1}sign(s_k-s_{k+1})\sum_{i=\alpha_k}^{\beta_k-1}\int_{x_i}^{x_{i+1}}\Big ( \hat{u}(x;k)-u_{ix}(k)(x-x_i)\Big )\Phi_i(k)dx =\nonumber\\ -\int_{0}^{T}\int_{0}^{s(t)}\tilde{u}^{\tau}\Phi_{t}^{\tau}\, dx \, dt 
-\int_{0}^{s_{1}}\phi^n(x)\Phi^{\tau}(x,\tau)\, dx -\int_{0}^{T-\tau}(\tilde{s}^{n}_1)^{\prime}(t)u^{\tau}((\tilde{s}^{n}_1)(t),t)\Phi^{\tau}((\tilde{s}^{n}_1)(t),t)\, dt\nonumber\\
-\sum_{k=1}^{n-1}\int_{t_{k-1}}^{t_{k}}\int_{s(t)}^{s_{k+1}}\tilde{u}^{\tau}\Phi_{t}^{\tau}\, dx\, dt
-\frac{1}{2}\sum_{k=1}^{n-1}sign(s_k-s_{k+1})\sum_{i=\alpha_k}^{\beta_k-1}h_i^2u_{ix}(k)\Phi_i(k),\label{Eq:W:3:32}
\end{gather}
where
\[ \alpha_k=\min(m_{j_k},m_{j_{k+1}}), \ \beta_k=\max(m_{j_k},m_{j_{k+1}}), \]
we derive that
\begin{gather}
\int_{0}^{T}\int_{0}^{s(t)}\bigg\{a \frac{\partial u^{\tau}}{\partial x}\Phi^{\tau}_x  - b \frac{\partial u^{\tau}}{\partial x}\Phi^{\tau}-c \tilde{u}^{\tau}\Phi^{\tau} + f \Phi^{\tau} -\tilde{u}^{\tau}\Phi_{t}^{\tau}\bigg\}\, dx \,dt-\int_{0}^{s_0}\phi^n(x)\Phi^{\tau}(x,\tau)\, dx\nonumber\\
-\int_{0}^{T-\tau}(\tilde{s}^{n}_1)^{\prime}(t)u^{\tau}((\tilde{s}^{n}_1)(t),t)\Phi^{\tau}((\tilde{s}^{n}_1)(t),t)\, dt
+\int_{0}^{T}g^n(t)\Phi^{\tau}(0,t)\, dt\nonumber\\
+\int_{0}^{T}\Big[\gamma(s^n(t),t)(s^n)^{\prime}(t)- \chi(s^n(t),t)) \Big]\Phi^{\tau}(s^n(t),t)\, dt
-R=0\label{Eq:W:3:32a}
\intertext{where}
R=\sum_{k=1}^{n}\int_{t_{k-1}}^{t_{k}}\int_{s(t)}^{s_{k}}\bigg\{a \frac{\partial u^{\tau}}{\partial x}\Phi^{\tau}_x  - b \frac{\partial u^{\tau}}{\partial x}\Phi^{\tau}-c \tilde{u}^{\tau}\Phi^{\tau} + f \Phi^{\tau}\bigg\}\, dx \,dt - \sum_{k=1}^{n-1}\int_{t_{k-1}}^{t_{k}}\int_{s(t)}^{s_{k+1}}\tilde{u}^{\tau}\Phi_{t}^{\tau}\, dx\, dt\nonumber\\
+\sum_{k=1}^{n}\int_{t_{k-1}}^{t_{k}}\int_{s^n(t)}^{s_{k}}\Big[ \gamma(s^n(t),t) (s^n)^{\prime}(t) - \chi(s^n(t),t))\Big]\frac{\partial \Phi^{\tau}}{\partial x}\, dx \, dt + \int_{s_0}^{s_1}\phi^n(x)\Phi^{\tau}(x,\tau)\, dx \nonumber\\+\sum_{k=1}^{n-1}\int_{t_{k-1}}^{t_{k}}\int_{s(t)}^{s_{k+1}}\tilde{u}^{\tau}\Phi_{t}^{\tau}\, dx\, dt
-\frac{1}{2}\sum_{k=1}^{n-1}sign(s_k-s_{k+1})\sum_{i=\alpha_k}^{\beta_k-1}h_i^2u_{ix}(k)\Phi_i(k)\nonumber
\end{gather}
First note that the sequence $\{\tilde{u}^\tau \}$ is equivalent to the sequence $\{u^\tau \}$ in strong, and accordingly also in a weak topology of $L_2(D)$, and hence converges to $u$ weakly in $L_2(D)$. Indeed, by using \eqref{firstenergyestimate1}
we have
\begin{gather}
\Vert  \tilde{u}^\tau - u^\tau \Vert^2_{L_2(D)} = \frac{1}{3}\sum_{k=1}^n \tau \sum_{i=0}^{N-1} h_i u_{ix}^2(k)  \max_{i} h_i^2  \to 0, \quad \text{as} \ n\to \infty. \label{Eq:W:3:32b} 
\end{gather}
Let
\[ \Delta = \bigcup_{k=1}^{n}\left\{ (x,t):t_{k-1}<t<t_{k},~\min(s(t),s_{k}) < x < \max(s(t),s_{k})\right\}\]
$|\Delta|$ denotes the Lebesgue measure of $\Delta$.
Since $\tilde{s}^n(t_k)=s_k$, we have
\begin{gather} 
|\Delta| \leq \sum_{k=1}^{n}\int_{t_{k-1}}^{t_{k}}\int_{t}^{t_{k}}|s'(\tau)|\, d\tau \, dt + \sum_{k=1}^n \tau |s(t_k)-\tilde{s}^n(t_k)| \leq \nonumber\\
\sqrt{T}\Vert s'\Vert_{L^{2}(0,T)}\tau + T\Vert s-\tilde{s}^n\Vert_{C[0,T]}\to 0 \quad \text{as}~\tau \to 0 \nonumber
\end{gather}
and all of the integrands are uniformly bounded in $L^{1}(D)$, it follows that the first term in the expression of $R$ converges to zero as $\tau \to 0$. In a similar way one can see that the second, third and fifth terms also converge to zero as $\tau \to 0$. The fourth term in the expression of $R$ converges to zero due to Corollary~\ref{lipshitz} and uniform convergence of $\{\Phi^\tau \}$ in $\overline{D}$. To prove the convergence to zero of the last term of $R$, let
\[ \tilde{\Delta} = \bigcup_{k=1}^{n-1}\left\{ (x,t):t_{k-1}<t<t_{k},d_{k}\equiv~\min(s_k,s_{k+1}) < x < d_{k+1}\equiv \max(s_k,s_{k+1})\right\}\]
From Corollary~\ref{lipshitz} it follows that
\[ |\tilde{\Delta}| \le C \tau \to 0, \quad \text{as}~\tau \to 0. \]
Since
\[ \sum_{i=\alpha_k}^{\beta_k-1}h_i = |s_k-s_{k+1}| \]
we have
\begin{gather}
\Big |\sum_{k=1}^{n-1}sign(s_k-s_{k+1})\sum_{i=\alpha_k}^{\beta_k-1}h_i^2u_{ix}(k)\Phi_i(k) \Big |
\le \sum_{k=1}^{n-1} |s_k-s_{k+1}|\int_{d_k}^{d_{k+1}} \Big | \frac{\partial u^\tau}{\partial x} \Big | |\Phi^\tau |dx \le \nonumber\\
C \sum_{k=1}^{n-1} \tau \int_{d_k}^{d_{k+1}} \Big | \frac{\partial u^\tau}{\partial x} \Big | |\Phi^\tau |dx
\le \Big \Vert \frac{\partial u^\tau}{\partial x} \Big \Vert_{L_2(\tilde{\Delta})}  \Vert \Phi^\tau  \Vert_{L_2(\tilde{\Delta})}\label{Eq:W:3:32c}
\end{gather}
Since the integrands are uniformly bounded in $L_2(D)$, the expression in \eqref{Eq:W:3:32c}
converges to zero as $\tau \to 0$.
Hence, we have
\begin{equation}
\lim_{\tau \to 0}R=0\label{Eq:W:3:33}
\end{equation}
Due to weak convergence of $u^\tau$ to $u$ in $W_2^{1,0}(D)$, weak convergence of $\tilde{u}^\tau$ to $u$ in $L_2(D)$ and uniform convergence of the sequences $\{ \Phi^{\tau}\}$, $\{ \frac{\partial \Phi^{\tau}}{\partial x}\}$ and $\{ \Phi_{t}^{\tau}\}$  to $\Phi$, $\frac{\partial \Phi}{\partial x}$ and $\frac{\partial \Phi}{\partial t}$ respectively, passing to limit as $\tau \to 0$, it follows that first, second and fourth integrals on the left-hand side of (\ref{Eq:W:3:32a}) converge to similar integrals with $u^\tau$ (or $\tilde{u}^\tau$), $\Phi^\tau$, $\Phi^\tau_t$, $\Phi^\tau(x,\tau)$, $g^n(t)$, $\phi^n(x)$ and $\Phi^\tau(0,t)$ replaced by $u$,$\Phi$, $\frac{\partial \Phi}{\partial t}$, $\Phi(x,0)$, $g(t)$, $\phi(x)$ and $\Phi(0,t)$ respectively. Since $s^n$ converges to $s$ strongly in $W_2^1[0,T]$, the traces $\gamma(s^n(t),(t))$, $\chi(s^n(t),t)$ converge strongly in $L_2[0,T]$ to traces $\gamma(s(t),(t))$, $\chi(s(t),t)$ respectively. Since $\Phi^{\tau}(s^n(t),t)$  converge uniformly on $[0,T]$ to
$\Phi(s(t),t)$, passing to the limit as $\tau \to 0$, the last integral on the left-hand side of 
(\ref{Eq:W:3:32a}) converge to similar integral with $s^n$ and $\Phi^\tau$ replaced by $s$ and $\Phi$.

It only remains to prove that
\begin{equation}
\lim_{\tau \to 0}\int_{0}^{T-\tau}(\tilde{s}^{n}_1)^{\prime}(t)u^{\tau}(\tilde{s}^{n}_1(t),t)\Phi^{\tau}(\tilde{s}^{n}_1(t),t)\, dt
=\int_{0}^{T}s'(t)u(s(t),t)\Phi(s(t),t)\label{Eq:W:3:34}
\end{equation}
Since $\{\tilde{s}^{n}_1\}$ converges to $s$ strongly in $W_2^1[0,T]$, from (\ref{Eq:W:3:27}) it follows that $\{ u^{\tau}(\tilde{s}^{n}_1(t),t)\}$ is uniformly bounded in $L_{2}[0,T]$ and
\begin{equation}\label{Eq:W:3:35}
\Vert u^{\tau}(\tilde{s}^{n}_1(t),t) - u^{\tau}(s(t),t)\Vert_{L_2[0,T]} \to 0\quad \text{as}~\tau \to 0
\end{equation}
Since $\{u^{\tau}\}$ converges to $u$ weakly in $W_2^{1,0}(D)$, it follows that
\begin{equation}
u^{\tau}(s(t),t) \to u(s(t),t),\quad \text{weakly in}~L_{2}[0,T]\label{Eq:W:3:35a}
\end{equation}
Since $\{\Phi^{\tau}(\tilde{s}^{n}_1(t),t)\}$ converges to $\Phi(s(t),t)$ uniformly in $[0,T]$, from
(\ref{Eq:W:3:35}),(\ref{Eq:W:3:35a}), (\ref{Eq:W:3:34}) easily follows.

Passing to the limit as $\tau \to 0$, from (\ref{Eq:W:3:32a}) it follows that $u$ satisfies integral identity (\ref{Eq:W:1:10}), i.e. it is a weak solution of the problem (\ref{Eq:W:1:1})-(\ref{Eq:W:1:4}). Since this solution is unique (\cite{LSU}) it follows that indeed the whole sequence $\{ u^{\tau}\}$ converges to $u \in V_2^{1,0}(\Omega)$ weakly in $W_2^{1,0}(\Omega)$. From the property of weak convergence and (\ref{Eq:W:3:30}),(\ref{V210estimate}) follows. Theorem is proved.

In particular, Theorem~\ref{continuity1} implies the following well-known existence result (\cite{LSU}):
\begin{corollary}\label{V210solution}
For arbitrary $v=(s,g)\in V_R$ there exists a weak solution $u \in V_2^{1,0}(\Omega)$ of the problem (\ref{Eq:W:1:1})-(\ref{Eq:W:1:4}) which satisfy the energy estimate (\ref{V210estimate})
\end{corollary}
\textbf{Remark}: All the proofs in this section can be pursued by using weaker assumptions
\[ \phi \in L_2[0,l], \ \gamma, \chi \in W_2^{1,0}(D), \ a \in L_{\infty}(D),  \]
and (\ref{Eq:W:1:6}) instead of conditions imposed in Section~\ref{E:1:4}. The only difference would be to define $\phi_i$ as a Steklov average
\[ \phi_i=\frac{1}{h_i}\int_{x_i}^{x_{i+1}}\phi(x)dx, \ i=0,...,N-1 \]
and replace the norm of $\phi^n$ in the first energy estimate through norm of $\phi$.
\subsection{Second Energy Estimate and its Consequences}\label{secondenergy}
Let given discrete control vector $[v]_n$, along with discrete state vector $[u([v]_n)]_n$, the vector 
\[ [\tilde{u}([v]_n)]_n=(\tilde{u}(0),\tilde{u}(1),...,\tilde{u}(n))  \]
is defined as
\[
\tilde{u}_i(k)=
\left\{
\begin{array}{l}
u_i(k) \ \ 0\le i \le m_j,\\
u_{m_j}(k) \ \ m_j < i \le N, k=\overline{0,n}.
\end{array}\right.
\]
The main goal of this section to prove the following energy estimation for the vector $\tilde{u}([v]_n)$.
\begin{theorem}\label{secondenergyestimate}
For all sufficiently small $\tau$ discrete state vector $[u([v]_n)]_n$ satisfies the following stability estimation:
\begin{gather}
\max_{1 \leq k \leq n}\sum_{i=0}^{m_j-1} h_i \tilde{u}_{ix}^2(k)+\tau \sum_{k=1}^{n}\sum_{i=0}^{m_j-1}h_i \tilde{u}_{i\overline{t}}^2(k)+\tau^2  \sum_{k=1}^{n}\sum_{i=0}^{m_j-1}h_i \tilde{u}_{ix\overline{t}}^2(k) \le \nonumber\\ C \bigg[ \left\Vert \phi^n\right\Vert_{L_{2}[0,s_0]}^{2}+\left\Vert \phi\right\Vert_{W_{2}^{1}[0,s_0]}^{2}+\left\Vert g^n\right\Vert_{W_{2}^{\frac{1}{4}}[0,T]}^{2}+\left\Vert \gamma(s^n(t),t)(s^n)'(t)\right\Vert_{W_{2}^{\frac{1}{4}}[0,T]}^{2}\nonumber\\ +\left\Vert \chi(s^n(t),t)\right\Vert_{W_{2}^{\frac{1}{4}}[0,T]}^{2}
+\left\Vert f\right\Vert_{L_{2}(D)}^{2}\bigg] ,\label{secondenergyestimate1} 
\end{gather}
\end{theorem}

\textbf{Proof}: Note that if $s_{k-1}\ge s_k$ then $u_i(k)$ can be replaced through $\tilde{u}_i(k)$ in all terms of (\ref{Eq:W:1:19}). By choosing $\eta_i=2\tau  \tilde{u}_{i\overline{t}}(k)$ in (\ref{Eq:W:1:19}) and by using the equality
\begin{gather}
2\tau a_{ik}\tilde{u}_{ix}(k)\tilde{u}_{ix\overline{t}}(k)=a_{ik}\tilde{u}_{ix}^2(k)-a_{i,k-1}\tilde{u}_{ix}^2(k-1)-\tau a_{ik\overline{t}} \tilde{u}_{ix}^2(k-1)
+\tau^2 a_{ik} \tilde{u}_{ik\overline{t}}^2(k) , \label{Eq:W:3:37}
\end{gather}
we have
\begin{gather}
\sum_{i=0}^{m_j-1}h_ia_{ik}\tilde{u}_{ix}^2(k)-\sum_{i=0}^{m_j-1}h_ia_{i,k-1}\tilde{u}_{ix}^2(k-1)+2\tau \sum_{i=0}^{m_j-1}h_i \tilde{u}_{i\overline{t}}^2(k)
+\tau^2 \sum_{i=0}^{m_j-1}h_i a_{ik} \tilde{u}_{ix\overline{t}}^2(k) \nonumber\\
=\tau  \sum_{i=0}^{m_j-1}h_i a_{ik\overline{t}} \tilde{u}_{ix}^2(k-1) + 2\tau \sum_{i=0}^{m_j-1}h_i b_{ik} \tilde{u}_{ix}(k) \tilde{u}_{i\overline{t}}(k) 
+2\tau \sum_{i=0}^{m_j-1}h_i c_{ik} \tilde{u}_{i}(k) \tilde{u}_{i\overline{t}}(k)\nonumber\\
-2\tau  \sum_{i=0}^{m_j-1}h_i f_{ik} \tilde{u}_{i\overline{t}}(k)-2\tau\left[(\gamma_{s^n} (s^n)')^{k}-\chi_{s^n}^{k} \right] \tilde{u}_{m_j,\overline{t}}(k)-2\tau g^n_k u_{0,\overline{t}}(k)\label{Eq:W:3:38}
\end{gather}
If $s_{k-1} < s_k$, then $u_i(k)$ can be replaced through $\tilde{u}_i(k)$ in all but in the term including backward discrete time derivative in (\ref{Eq:W:1:19}). The latter will be estimated with the help of the following inequality:
\begin{gather}
2\tau\sum_{i=0}^{m_j-1}h_iu_{i\overline{t}}(k)\tilde{u}_{i\overline{t}}(k)\ge \tau\sum_{i=0}^{m_j-1}h_i\tilde{u}_{i\overline{t}}^2(k)-(C^\prime)^2 \tau \sum_{i=m_{j_{k-1}}}^{m_j-1}h_i u_{ix}^2(k-1)\label{Eq:W:3:38a}
\end{gather}
To prove \eqref{Eq:W:3:38a}, we transform the left hand side with the help of the CBS and Cauchy inequalty with $\epsilon=\tau$ to derive
\begin{gather}
2\tau\sum_{i=0}^{m_j-1}h_iu_{i\overline{t}}(k)\tilde{u}_{i\overline{t}}(k)=2\tau\sum_{i=0}^{m_j-1}h_i\tilde{u}_{i\overline{t}}^2(k)-2\sum_{i=m_{j_{k-1}}}^{m_j-1}h_i\tilde{u}_{i\overline{t}}(k)\sum_{p=m_{j_{k-1}}}^{i-1}h_p u_{px}(k-1)\nonumber\\
\ge \tau\sum_{i=0}^{m_j-1}h_i\tilde{u}_{i\overline{t}}^2(k) - \frac{1}{\tau}\sum_{i=m_{j_{k-1}}}^{m_j-1}h_i \Big ( \sum_{p=m_{j_{k-1}}}^{i-1}h_p u_{px}(k-1)\Big )^2 \nonumber\\
\ge \tau\sum_{i=0}^{m_j-1}h_i\tilde{u}_{i\overline{t}}^2(k) - \frac{1}{\tau} |s_k-s_{k-1}|^2\sum_{i=m_{j_{k-1}}}^{m_j-1}h_i u_{ix}^2(k-1) \label{Eq:W:3:38b}
\end{gather}
which implies \eqref{Eq:W:3:38a} due to \eqref{Eq:W:2:28}. Hence, in general \eqref{Eq:W:3:38} is replaced with the inequality
\begin{gather}
\sum_{i=0}^{m_j-1}h_ia_{ik}\tilde{u}_{ix}^2(k)-\sum_{i=0}^{m_j-1}h_ia_{i,k-1}\tilde{u}_{ix}^2(k-1)+\tau \sum_{i=0}^{m_j-1}h_i \tilde{u}_{i\overline{t}}^2(k)
+\tau^2 \sum_{i=0}^{m_j-1}h_i a_{ik} \tilde{u}_{ix\overline{t}}^2(k) \nonumber\\
\le \tau \sum_{i=0}^{m_j-1}h_i a_{ik\overline{t}} \tilde{u}_{ix}^2(k-1) + 2\tau \sum_{i=0}^{m_j-1}h_i b_{ik} \tilde{u}_{ix}(k) \tilde{u}_{i\overline{t}}(k) 
+2\tau \sum_{i=0}^{m_j-1}h_i c_{ik} \tilde{u}_{i}(k) \tilde{u}_{i\overline{t}}(k) \nonumber\\ +(C^\prime)^2 \tau\sum_{i=m_{j_{k-1}}}^{m_j-1}h_i u_{ix}^2(k-1)
-2\tau  \sum_{i=0}^{m_j-1}h_i f_{ik} \tilde{u}_{i\overline{t}}(k)\nonumber\\ -2\tau\left[(\gamma_{s^n} (s^n)')^{k}-\chi_{s^n}^{k} \right] \tilde{u}_{m_j,\overline{t}}(k)-2\tau g^n_k \tilde{u}_{0,\overline{t}}(k)\label{Eq:W:3:38c}
\end{gather}
By adding inequalities (\ref{Eq:W:3:38c}) with respect to $k$ from 1 to arbitrary $p \leq n$ we derive
\begin{gather}
\sum_{i=0}^{m_{j_p}-1}h_ia_{ip}\tilde{u}_{ix}^2(p)+\tau\sum_{k=1}^{p}  \sum_{i=0}^{m_j-1}h_i \tilde{u}_{i\overline{t}}^2(k)+\tau^2\sum_{k=1}^p  \sum_{i=0}^{m_j-1}h_i a_{ik} \tilde{u}_{ix\overline{t}}^2(k)\nonumber\\
\leq (C^\prime)^2 \tau\sum_{k=1}^p 1_+(s_k-s_{k-1})\sum_{i=m_{j_{k-1}}}^{m_j-1}h_i u_{ix}^2(k-1)+\tau \sum_{k=1}^p\sum_{i=0}^{m_j-1}h_i a_{ik\overline{t}} \tilde{u}_{ix}^2(k-1)  \nonumber\\   + 2\tau\sum_{k=1}^p \Big [ \sum_{i=0}^{m_j-1}h_i b_{ik} \tilde{u}_{ix}(k) \tilde{u}_{i\overline{t}}(k) 
+2 \sum_{i=0}^{m_j-1}h_i c_{ik} \tilde{u}_{i}(k) \tilde{u}_{i\overline{t}}(k) - 2 \sum_{i=0}^{m_j-1}h_i f_{ik} \tilde{u}_{i\overline{t}}(k) \Big ]
\nonumber\\
+\sum_{i=0}^{m_{j_0}-1}h_ia_{i0}\phi_{ix}^2-2\tau\sum_{k=1}^{p}\left[ (\gamma_{s^n} (s^n)')^{k}-\chi_{s^n}^{k}\right] \tilde{u}_{m_j,\overline{t}}(k)-2\tau \sum_{k=1}^{p} g^n_k \tilde{u}_{0,\overline{t}}(k)\label{Eq:W:3:39}
\end{gather}
By using (\ref{Eq:W:1:6}) and by applying Cauchy inequalities with appropriately chosen
$\ep > 0$, from (\ref{Eq:W:3:39}) it follows that
\begin{gather}
a_0\sum_{i=0}^{m_{j_p}-1}h_i\tilde{u}_{ix}^2(p)+\frac{\tau}{2}\sum_{k=1}^{p}  \sum_{i=0}^{m_j-1}h_i \tilde{u}_{i\overline{t}}^2(k)+a_0\tau^2\sum_{k=1}^p  \sum_{i=0}^{m_j-1}h_i  \tilde{u}_{ix\overline{t}}^2(k)\le \tau \sum_{k=1}^p\sum_{i=0}^{m_j-1}h_i a_{ik\overline{t}} \tilde{u}_{ix}^2(k-1)\nonumber\\
+ C \tau \sum_{k=1}^{n}\bigg[\sum_{i=0}^{m_j-1}h_iu_i^{2}(k) + \sum_{i=0}^{m_j-1}h_iu_{ix}^2(k) + \sum_{i=0}^{m_j-1}h_i f_{ik}^2 \bigg] + C \sum_{i=0}^{m_{j_0}-1}h_i\phi_{ix}^2 \nonumber\\ + 2\tau \sum_{k=1}^{n}\left|(\gamma_{s^n} (s^n)')^{k}-\chi_{s^n}^{k} \right|\left|\tilde{u}_{m_j,\overline{t}}(k)\right| + 2\tau\sum_{k=1}^{n}\left|g^n_{k}\right| \left|\tilde{u}_{0,\overline{t}}(k)\right|\label{Eq:W:3:40}
\end{gather}
where $C$ is independent of $n$. First term on the right hand side will be estimated as follows:
\begin{gather}
\tau \sum_{k=1}^p\sum_{i=0}^{m_j-1}h_i a_{ik\overline{t}} \tilde{u}_{ix}^2(k-1) = \sum_{k=1}^n \sum_{i=0}^{m_j-1} \frac{1}{\tau}\int_{x_i}^{x_{i+1}}\int_{t_{k-1}}^{t_k}\int_{t-\tau}^{t} \frac{\partial a(x,\xi)}{\partial \xi} d\xi dt dx u_{ix}^2(k-1) \nonumber\\
\le 2  \int_{0}^{T}esssup_{0 \leq x \leq l} \left| \frac{\partial a(x,t)}{\partial t}\right| dt  \max_{1 \leq k \leq n}\sum_{i=0}^{m_j-1} h_i \tilde{u}_{ix}^2(k) + C \sum_{i=0}^{m_{j_0}-1}h_i\phi_{ix}^2
\label{Eq:W:3:40a}
\end{gather}
Due to arbitraricity of $p$, from \eqref{Eq:W:3:40} it follows
\begin{gather}
a_0 \max_{1\le k \le n}\sum_{i=0}^{m_{j}-1}h_i\tilde{u}_{ix}^2(p)+\frac{\tau}{2}\sum_{k=1}^{n}  \sum_{i=0}^{m_j-1}h_i \tilde{u}_{i\overline{t}}^2(k)+a_0\tau^2\sum_{k=1}^n  \sum_{i=0}^{m_j-1}h_i  \tilde{u}_{ix\overline{t}}^2(k)\nonumber\\ \le 
 2  \int_{0}^{T}esssup_{0 \leq x \leq l} \left| \frac{\partial a(x,t)}{\partial t}\right| dt  \max_{1 \leq k \leq n}\sum_{i=0}^{m_j-1} h_i \tilde{u}_{ix}^2(k) + C \sum_{i=0}^{m_{j_0}-1}h_i\phi_{ix}^2
\nonumber\\
+  C \tau \sum_{k=1}^{n}\bigg[\sum_{i=0}^{m_j-1}h_iu_i^{2}(k) + \sum_{i=0}^{m_j-1}h_iu_{ix}^2(k) + \sum_{i=0}^{m_j-1}h_i f_{ik}^2 \bigg]  \nonumber\\ + 2\tau \sum_{k=1}^{n}\left|(\gamma_{s^n} (s^n)')^{k}-\chi_{s^n}^{k} \right|\left|\tilde{u}_{m_j,\overline{t}}(k)\right| + 2\tau\sum_{k=1}^{n}\left|g^n_{k}\right| \left|\tilde{u}_{0,\overline{t}}(k)\right|\label{Eq:W:3:40b}
\end{gather}
If
\begin{gather}
 2  \int_{0}^{T}esssup_{0 \leq x \leq l} \left| \frac{\partial a(x,t)}{\partial t}\right| dt < a_0
\label{Eq:W:3:40c}
\end{gather}
then the first term on the right hand side of \eqref{Eq:W:3:40b} is absorbed into the first term on the left hand side. If \eqref{Eq:W:3:40c} is not satisfied, then we can partition $[0,T]$ into finitely many subsegments which obey \eqref{Eq:W:3:40c}, absorb first term on the right hand side into the left hand side in each subsegment and through summation achieve the same for  \eqref{Eq:W:3:40b} in general. Hence we have
\begin{gather}
\max_{1\le k \le n}\sum_{i=0}^{m_{j}-1}h_i\tilde{u}_{ix}^2(p)+\tau\sum_{k=1}^{n}  \sum_{i=0}^{m_j-1}h_i \tilde{u}_{i\overline{t}}^2(k)+\tau^2\sum_{k=1}^n  \sum_{i=0}^{m_j-1}h_i  \tilde{u}_{ix\overline{t}}^2(k)\nonumber\\ \le 
 C \tau \sum_{k=1}^{n}\bigg[\sum_{i=0}^{m_j-1}h_iu_i^{2}(k) + \sum_{i=0}^{m_j-1}h_iu_{ix}^2(k) + \sum_{i=0}^{m_j-1}h_i f_{ik}^2 \bigg] + C \sum_{i=0}^{m_{j_0}-1}h_i\phi_{ix}^2 \nonumber\\ + C\tau \sum_{k=1}^{n}\left|(\gamma_{s^n} (s^n)')^{k}-\chi_{s^n}^{k} \right|\left|\tilde{u}_{m_j,\overline{t}}(k)\right| + C\tau\sum_{k=1}^{n}\left|g^n_{k}\right| \left|\tilde{u}_{0,\overline{t}}(k)\right|\label{Eq:W:3:40d}
\end{gather}
with some $C$ independent of $n$.

Since $\gamma, \chi \in W_2^{1,1}(D)$ we have $\gamma(s^n(t),t), \chi(s^n(t),t) \in W_2^{\frac{1}{4}}[0,T]$ (\cite{Nikolski, BIN, LSU}) and
\begin{equation}\label{Eq:W:3:41}
\Vert \gamma(s^n(t),t)\Vert_{W_2^{\frac{1}{4}}[0,T]}\leq C \Vert \gamma \Vert_{W_2^{1,1}(D)}, \  \Vert \chi(s^n(t),t)\Vert_{W_2^{\frac{1}{4}}[0,T]}  \leq C \Vert \chi \Vert_{W_2^{1,1}(D)},
\end{equation}
where $C$ is independent of $n$. According to Lemma~\ref{mappings} ${\mathcal P}_n([v]_n) \in V_{R+1}$. By applying Morrey inequality to $(s^n)^{\prime}$ we easily deduce that $\gamma(s^n(t),t)(s^n)^{\prime}(t) \in W_2^{\frac{1}{4}}[0,T]$ and moreover,
\begin{equation}\label{Eq:W:3:42}
\Vert \gamma(s^n(t),t)(s^n)^\prime(t)\Vert_{W_2^{\frac{1}{4}}[0,T]}\leq C_1 \Vert \gamma(s^n(t),t)\Vert_{W_2^{\frac{1}{4}}[0,T]} \Vert s^n\Vert_{W_2^2[0,T]} \leq C \Vert \gamma \Vert_{W_2^{1,1}(D)},
\end{equation}
where $C$ is independent of $n$.

Let $w(x,t)$ be a function in $W_{2}^{2,1}(D)$ such that
\begin{equation}\label{Eq:W:3:43}
w(x,0)=\phi(x) \quad \text{for}~x\in [0,s_0], \ a(0,t)w_{x}(0,t)=g^n(t),\ \quad \text{for a.e.}~t\in [0,T]
\end{equation}
\begin{equation}\label{Eq:W:3:44}
a(s^n(t),t)w_{x}(s^n(t),t)=\gamma(s^n(t),t)(s^n)^{\prime}(t)-\chi(s^n(t),t) \quad \text{for a.e.}~t\in [0,T]
\end{equation}
and
\begin{gather}
\left\Vert w\right\Vert_{W_{2}^{2,1}(D)} \leq C\Big [ \left\Vert g^n\right\Vert_{W_{2}^{\frac{1}{4}}[0,T]} +\left\Vert \phi(x)\right\Vert_{W_{2}^{1}[0,s_{0}]}\nonumber\\
+\left\Vert \gamma(s^n(t),t)(s^n)^{\prime}(t)-\chi(s^n(t),t)\right\Vert_{W_{2}^{\frac{1}{4}}[0,T]}\Big ]\label{Eq:W:3:45}
\end{gather}
The existence of $w$ follows from the result on traces of Sobolev functions \cite{BIN,Nikolski}.
For example, $w$ can be constructed as a solution from $W_2^{2,1}(\Omega^n)$ of the heat equation 
in 
\[ \Omega^n=\{0<x<s^n(t), 0<t<T \} \]
under initial-boundary conditions (\ref{Eq:W:3:43}),(\ref{Eq:W:3:44})with subsequent continuation to $W_2^{2,1}(D)$ with norm preservation \cite{Solonnikov1, Solonnikov2}. 

Hence, by replacing in the original problem (\ref{Eq:W:1:1})-(\ref{Eq:W:1:4}) $u$ with $u-w$ we can derive modified (\ref{Eq:W:3:40d}) without the last three terms on the right-hand side and with $f$, replaced by
\begin{equation}\label{Eq:W:3:46}
F=f+w_{t}-(a w_x)_x-bw_x-cw\in L_{2}(D).
\end{equation}
By using the stability estimation (\ref{firstenergyestimate3}),  from modified (\ref{Eq:W:3:40d}),(\ref{Eq:W:3:45}) and (\ref{Eq:W:3:46}),
the following estimation follows:
\begin{gather}
\max_{1 \leq k \leq n}\sum_{i=0}^{m_j-1} h_i \tilde{u}_{ix}^2(k)+\tau \sum_{k=1}^{n}\sum_{i=0}^{m_j-1}h_i \tilde{u}_{i\overline{t}}^2(k)+\tau^2  \sum_{k=1}^{n}\sum_{i=0}^{m_j-1}h_i \tilde{u}_{ix\overline{t}}^2(k)   \le \nonumber\\ C \bigg[\left\Vert \phi^n\right\Vert_{L_{2}[0,s_0]}^{2} + \left\Vert \phi\right\Vert_{W_{2}^{1}[0,s_0]}^{2}+\left\Vert g^n\right\Vert_{W_{2}^{\frac{1}{4}}[0,T]}^{2}+\left\Vert \gamma(s^n(t),t)(s^n)'(t)\right\Vert_{W_{2}^{\frac{1}{4}}[0,T]}^{2}+\nonumber\\ \left\Vert \chi(s^n(t),t)\right\Vert_{W_{2}^{\frac{1}{4}}[0,T]}^{2}
+\left\Vert f\right\Vert_{L_{2}(D)}^{2} +   \sum_{k=1}^{n-1}{\bf 1}_{+}(s_{k+1}-s_{k}) \sum_{i=m_j}^{m_{j_{k+1}}-1} h_iu_i^2(k)  \bigg] ,\label{secondenergyestimate2} 
\end{gather}
By estimating the last term on the right hand side of \eqref{secondenergyestimate2} as in the proof of Theorem~\ref{continuity1}, 
(\ref{secondenergyestimate1}) follows.  Theorem is proved.

Second energy estimate (\ref{secondenergyestimate1}) allows to strengthen the result of Theorem~\ref{continuity1}. 
\begin{theorem}\label{continuity2}
Let $[v]_n \in V_R^n, n=1,2,...$ be a sequence of discrete controls and the sequence $\{\mathcal{P}_n([v]_n\}$ converges weakly in $W_2^2[0,T] \times W_2^1[0,T]$ to $v=(s,g)$ (i.e. strongly in $W_2^1[0,T] \times L_2[0,T]$) to $v=(s,g)$ for any $\delta>0$,
  \[
    \Omega^{\prime} = \Omega \cap \{x<s(t)-\delta,~0<t<T\}
  \] 
Then the sequence $\{\hat{u}^\tau\}$ converges as $\tau \to 0$ weakly in $W_2^{1,1}(\Omega^{\prime})$ to weak solution $u \in W_2^{1,1}(\Omega)$ of the problem (\ref{Eq:W:1:1})-(\ref{Eq:W:1:4}), i.e. to the solution of the integral identity (\ref{Eq:W:1:9}). Moreover, $u$ satisfies the energy estimate
\begin{equation}\label{W211estimate}
\Vert u \Vert_{W_2^{1,1}(\Omega)}^2 \le C \Big ( \Vert \phi \Vert_{W_2^1(0,s_0)}^2 + \Vert g \Vert_{W_2^{\frac{1}{4}}[0,T]}^2 + \Vert f \Vert_{L_2(D)}^2 + 
\Vert \gamma \Vert_{W_2^{1,1}(D)}^2 + \Vert \chi \Vert_{W_2^{1,1}(D)}^2 \Big )
\end{equation}
\end{theorem}

\textbf{Proof}: Let $\epsilon_m \downarrow 0$ be arbitrary sequence and 
\[ \Omega_m=\{(x,t): 0<x<s(t)-\epsilon_m, 0<t\le T\}  \]
Note that the sequence $s^n$ converges to $s$ uniformly in $[0,T]$. 
  Calculate
 \[   \Vert \hat{u}^\tau \Vert_{W_2^{1,1}(\Omega_m)}^2
       = \sum_{k = 1}^n \int_{t_{k-1}}^{t_k} \int_0^{s(t) - \epsilon_m} |\hat{u}^\tau|^2 + \Big | \frac{\partial{\hat{u}^\tau}}{\partial x}\Big |^2 + \Big | \frac{\partial{\hat{u}^\tau}}{\partial t}\Big |^2 \,dx \,dt \]
Denote $s^m_k = x_{\hat{i}}$, where
\[    \hat{i}=\max \Big \{ i\leq N: \max_{t_{k-1} \leq t \leq t_{k}} s(t) - \epsilon_m\leq x_i \leq \max_{t_{k-1} \leq t \leq t_{k}} s(t) -  \frac{\epsilon_m}{2} \Big \} \]
and substitute to derive
  \begin{align}
   \Vert \hat{u}^\tau \Vert_{W_2^{1,1}(\Omega_m)}^2
      & \leq \sum_{k = 1}^n \int_{t_{k-1}}^{t_k} \int_0^{s^m_k} | \hat{u}(x;k-1)+ \hat{u}_{\bar{t}}(x;k)(t - t_{k-1}) |^2 \,dx \,dt + \nonumber
    \\
      & \quad + \sum_{k = 1}^n \int_{t_{k-1}}^{t_k} \int_0^{s^m_k} \Big | \frac{d \hat{u}(x;k-1)}{dx}+\frac{d \hat{u}_{\bar{t}}(x;k)}{dx}(t - t_{k-1}) \Big |^2 \,dx \,dt +\nonumber
    \\
      & \quad + \sum_{k = 1}^n \int_{t_{k-1}}^{t_k} \int_0^{s^m_k} | \hat{u}_{\bar{t}}(x;k)|^2 \,dx \,dt \nonumber
    \intertext{From which by using \eqref{htau} it follows that}
     \Vert \hat{u}^\tau \Vert_{W_2^{1,1}(\Omega_m)}^2
      & \leq C\Big \{ \tau \sum_{k = 1}^n \sum_0^{\hat{i}-1} h_i u_i^2(k-1) +  \tau \sum_{k = 1}^n \sum_0^{\hat{i}-1} h_i u_{ix}^2(k-1) + \nonumber
    \\
      & \quad +  \tau \sum_{k = 1}^n \sum_0^{\hat{i}-1} h_i u_{i{\bar{t}}}^2(k)+\tau^2 \sum_{k = 1}^n \sum_0^{\hat{i}-1} h_i u_{ix{\bar{t}}}^2(k)\Big \}\label{357}
  \end{align}
  where $C$ is independent of $n$ or $\tau$. Our goal is to prove that the right hand side of \eqref{357} is bounded by the left hand side of \eqref{secondenergyestimate1} for sufficiently large $n$. It is sufficient to prove the following claim: {\it 
    for fixed $\epsilon_m$, there exists $N=N(\epsilon_m)$ such that for $\forall n> N$}
 \begin{equation}
      s_k^m < \min(s_k, s_{k-1}),~k=1,\ldots,n\label{omegaminterior}
    \end{equation}
Indeed from \eqref{omegaminterior} it follows that 
      \begin{align}
      \tau \sum_{k = 1}^n \sum_0^{\hat{i}-1} h_i u_{ix}^2(k-1)  +  \tau \sum_{k = 1}^n \sum_0^{\hat{i}-1} h_i u_{i{\bar{t}}}^2(k)+\tau^2 \sum_{k = 1}^n \sum_0^{\hat{i}-1} h_i u_{ix{\bar{t}}}^2(k)\nonumber\\
    \leq   \tau \sum_{k = 1}^n \sum_0^{m_j-1} h_i \tilde{u}_{ix}^2(k-1)  +  \tau \sum_{k = 1}^n \sum_0^{m_j-1} h_i \tilde{u}_{i{\bar{t}}}^2(k)+\tau^2 \sum_{k = 1}^n \sum_0^{m_j-1} h_i \tilde{u}_{ix{\bar{t}}}^2(k)\label{compactness}
      \end{align}
To prove \eqref{omegaminterior}, we first show that for sufficiently large $n$ and all $t_{k-1} \leq t \leq t_{k}$
    \begin{equation}
      s(t) - s_k < \frac{\epsilon_m}{2}\label{eq:second-energy-est-conseq-smk-term-enough-1}
    \end{equation}
    We have
    \begin{gather*}
      s(t) - s_k
      = s(t) - s(t_k) + s(t_k) - s^n(t_k) + s^n(t_k) - s_k\nonumber
      \\
      \leq \Vert s' \Vert_{C[0,T]} \tau
      + \Vert s-s^n \Vert_{C[0,T]}
      + s^n(t_k) - s_k\nonumber
    \end{gather*}
    Observe that
    \[
      s^n(t_k) = \frac{s_k + s_{k-1}}{2}\nonumber
    \]
    so
    \[|s(t) - s_k|
      \leq \Vert s' \Vert_{C[0,T]} \tau
      + \Vert s-s^n \Vert_{C[0,T]}
      + \frac{|s_{k-1} - s_k|}{2}\nonumber
    \]
    By Morrey's inequality and~\eqref{lipshitz},
    \begin{equation}\label{alma}
      |s(t) - s_k|
      \leq \left( C\Vert s' \Vert_{W_2^1[0,T]} + \frac{C'}{2}\right)\tau
      + \Vert s-s^n \Vert_{C[0,T]}
    \end{equation}
    By uniform convergence of $s^n \to s$ on $[0,T]$ there exists $N_1=N_1(\epsilon_m)$ for which~\eqref{eq:second-energy-est-conseq-smk-term-enough-1} holds.
    Similarly, there exists some $N_2=N_2(\epsilon_m)$, such that for $n>N_2$ and for all $t_{k-1}\leq t\leq t_{k}$, \eqref{eq:second-energy-est-conseq-smk-term-enough-1} is true with $s_k$ replaced by $s_{k-1}$. Hence, \eqref{omegaminterior} follows with $N=\max(N_1,N_2)$.
 
  Applying~\eqref{357}, \eqref{compactness}, second energy estimate~\eqref{secondenergyestimate1},
  and the first energy estimate~\eqref{firstenergyestimate1}, derive
  \begin{gather}
    \Vert \hat{u}^\tau \Vert_{W_2^{1,1}(\Omega_m)}^2
    \leq
    C\Big[ \left\Vert \phi^n\right\Vert_{L_{2}[0,s_0]}^{2} +
    \Vert \phi \Vert_{W_{2}^{1}[0,s_{0}]}^{2}
    + \Vert f \Vert_{L_{2}(D)}^{2}
    + \Vert \gamma \Vert_{W_{2}^{1,1}(D)}^{2} + \nonumber
    \\
    + \Vert \chi \Vert_{W_{2}^{1,1}(D)}^{2}
    + \Vert g^n \Vert_{W_{2}^{1/4}[0,T]}^{2}
    + \sum_{k=1}^{n-1}{\bf 1}_{+}(s_{k+1}-s_{k}) \sum_{i=m_j}^{m_{j_{k+1}}-1} h_iu_i^2(k)
    \Big]\label{eq:second-energy-est-conseq-2}
  \end{gather}
    By estimating the last term on the right hand side of \eqref{eq:second-energy-est-conseq-2} as in the proof of Theorem~\ref{continuity1}, we derive
\begin{gather}\label{Eq:W:3:51}
\Vert \hat{u}^\tau \Vert_{W_2^{1,1}(\Omega_m)}^2 \le C \Big ( \Vert \phi^n \Vert_{L_2(0,s_0)}^2+\Vert \phi \Vert_{W_2^1(0,s_0)}^2 + \Vert g^n \Vert_{W_2^{\frac{1}{4}}[0,T]}^2 + \Vert f \Vert_{L_2(D)}^2 + \nonumber\\
\Vert \gamma \Vert_{W_2^{1,1}(D)}^2 + \Vert \chi \Vert_{W_2^{1,1}(D)}^2 \Big ), \label{Eq:W:3:51}
\end{gather}
Since $\phi^n \rightarrow \phi$ strongly in $L_2[0,s_0]$, $g^n \rightarrow g$ weakly in $W_2^1[0,T]$, the right hand side is uniformly bounded independent of $n$. 
Hence, $\{\hat{u}^\tau\}$ is weakly precompact in $W_2^{1,1}(\Omega_m)$. It follows that it is strongly precompact
in $L_2(\Omega_m)$. Let $u$ be a weak limit point of $\{\hat{u}^\tau\}$ in $W_2^{1,1}(\Omega_m)$, and therefore a strong limit point in $L_2(\Omega_m)$. From anothe side the sequences $\{\hat{u}^\tau\}$ and $\{u^\tau \}$ are equivalent in strong topology of $L_2(\Omega_m)$. Indeed, we have for all $n>N(m)$
\begin{equation}\label{armud} 
\Vert \hat{u}^\tau - u^\tau \Vert_{L_2(\Omega_m)}^2\le 2\tau^3 \sum_{k=1}^{n} \sum_{i=0}^{m_j-1} \Big [ h_i  \tilde{u}_{i\bar{t}}^2(k) +\frac{1}{3}h_i^3 \tilde{u}_{ix\overline{t}}^2(k) \Big ] =O(\tau), \quad \text{as}~\tau \to 0, 
\end{equation}
due to second energy estimate (\ref{secondenergyestimate1}).
Therefore, $u$ is a strong limit point of the sequence $\{u^\tau\}$ in $L_2(\Omega_m)$. By Theorem~\ref{continuity1} whole sequence $\{u^\tau\}$
converges weakly in $W_2^{1,0}(\Omega)$ to the unique weak solution from $V_2^{1,0}(\Omega)$ of the problem (\ref{Eq:W:1:1})-(\ref{Eq:W:1:4}). Hence, $u$ is a weak solution of the problem (\ref{Eq:W:1:1})-(\ref{Eq:W:1:4}) and we conclude that whole sequence $\{\hat{u}^\tau\}$ converges weakly in $W_2^{1,1}(\Omega_m)$ to $u \in W_2^{1,1}(\Omega_m)$ which is a weak solution of the problem (\ref{Eq:W:1:1})-(\ref{Eq:W:1:4}) from $V_2^{1,0}(\Omega)$. Hence, $u_t$ exists in $\Omega_m$ and $\Vert u_t\Vert_{L_2(\Omega_m)}$ is uniformly bounded by the right hand side of \eqref{Eq:W:3:51}. It easily follows that the weak derivative $u_t$ exists in $\Omega$, and $u \in W_2^{1,1}(\Omega)$. By using well-known property of the weak convergence, and passing to limit first as $n \rightarrow +\infty$, and then as $m \rightarrow +\infty$, from \eqref{Eq:W:3:51}, \eqref{W211estimate} follows. Theorem is proved.

In particular, Theorem~\ref{continuity2} implies the following  existence result:
\begin{corollary}\label{W211solution}
For arbitrary $v=(s,g)\in V_R$ there exists a weak solution $u \in W_2^{1,1}(\Omega)$ of the problem (\ref{Eq:W:1:1})-(\ref{Eq:W:1:4}) which satisfy the energy estimate (\ref{W211estimate}). By Sobolev extension theorem $u$ can be continued to $W_2^{1,1}(D)$ with the norm preservation:
\begin{equation}\label{W211estimate-1}
\Vert u \Vert_{W_2^{1,1}(D)}^2 \le C \Big ( \Vert \phi \Vert_{W_2^1(0,s_0)}^2 + \Vert g \Vert_{W_2^{\frac{1}{4}}[0,T]}^2 + \Vert f \Vert_{L_2(D)}^2 + 
\Vert \gamma \Vert_{W_2^{1,1}(D)}^2 + \Vert \chi \Vert_{W_2^{1,1}(D)}^2 \Big )
\end{equation}
\end{corollary}
\textbf{Remark}: In fact, we proved slightly higher regularity of $u$, and both in Theorem~\ref{continuity2} 
and Corollary~\ref{W211solution} $W_2^{1,1}(\Omega)$ or $W_2^{1,1}(D)$-norm on the left-hand sides of (\ref{W211estimate}) or \eqref{W211estimate-1} can be replaced with 
\[ {\bf \Vert} u {\bf \Vert}^2 = \max_{0\le t \le T}\Vert u(x,t)\Vert_{W_2^1[0,s(t)]}^2+\Vert u_t \Vert_{L_2(\Omega)}^2 \quad \text{or}~  {\bf \Vert} u {\bf \Vert}^2 = \max_{0\le t \le T}\Vert u(x,t)\Vert_{W_2^1[0,l]}^2+\Vert u_t \Vert_{L_2(D)}^2 \]

The proof of the Theorem~\ref{existence} coincides with the proof of identical Theorem~\ref{existence} in \cite{Abdulla1}. The main idea is that  first and second energy estimates imply weak continuity of the functional $\mathcal{J}(v)$ in $W_2^2[0,T] \times W_2^1[0,T]$. Since $V_R$ is weakly compact existence of the optimal control follows from Weierstrass theorem in weak topology.

We split the remainder of the proof of Theorem~\ref{convergence} into three lemmas.
\begin{lemma}\label{condition(3)}
Let $\mathcal{J}_{*}(\pm \epsilon) = \inf\limits_{V_{R\pm \epsilon}} \mathcal{J}(v)$, $\ep >0$. Then
\begin{equation}
\lim_{\epsilon \to 0} \mathcal{J}_{*}(\epsilon) = \mathcal{J}_{*} = \lim_{\epsilon \to 0}\mathcal{J}_{*}(-\epsilon)\label{Eq:W:3:55}
\end{equation}
\end{lemma}
The proof of Lemma~\ref{condition(3)} coincides with the proof of identical lemma 3.9 from \cite{Abdulla1}.
\begin{lemma}\label{condition(1)}
For arbitrary $v=(s,g) \in V_{R}$,
\begin{equation}\label{Eq:W:3:58}
\lim\limits_{n\to \infty} \mathcal{I}_n(\mathcal{Q}_n(v))=\mathcal{J}(v)
\end{equation}
\end{lemma}
\textbf{Proof:}
Let $v \in V_{R}$, $u=u(x,t;v)$, $\mathcal{Q}_{n}(v)=[v]_{n}$ and $[u([v]_{n})]_{n}$ be a corresponding discrete state vector. In Theorem~\ref{continuity2} it is proved that the sequence $\{\hat{u}^\tau\}$ converges to $u$ weakly in $W_2^{1,1}(\Omega_m)$ for any fixed $m$. This implies that the sequences of traces $\{\hat{u}^\tau(0,t)\}$ and $\{\hat{u}^\tau(s(t)-\epsilon_m,t)\}$ converge strongly in $L_2[0,T]$ to corresponding traces $u(0,t)$ and $u(s(t)-\epsilon_m,t)$.
Let us prove that that the sequences of traces $\{u^{\tau}(0,t)\}$ and $\{u^{\tau}(s(t)-\epsilon_m,t)\}$ converge strongly in $L^{2}[0,T]$ to traces $u(0,t)$ and $u(s(t)-\epsilon_m,t)$ respectively. By Sobolev embedding theorem (\cite{BIN,Nikolski}) it is enough to prove that the sequences $\{u^\tau\}$ and $\{\hat{u}^\tau\}$ are equivalent in strong topology of $W_2^{1,0}(\Omega_m)$. In Theorem~\ref{continuity2} it is proved that they are equivalent in strong topology of $L_2(\Omega_m)$. It remains only to demonstrate that the sequences of 
derivatives $\frac{\partial u^\tau}{\partial x}$ and $\frac{\partial\hat{u}^\tau}{\partial x}$ are equivalent in strong topology of $L_2(\Omega_m)$. 
Following the proof of the Theorem~\ref{continuity2}, from the second energy estimate (\ref{secondenergyestimate1}) it follows that for all $n>N(m)$ 
\begin{equation}\label{Eq:W:3:59}
\Big \Vert \frac{\partial u^\tau}{\partial x}-\frac{\partial\hat{u}^\tau}{\partial x}\Big \Vert_{L_2(\Omega_m)}^2\le \frac{1}{3}\sum_{k=1}^n \tau^3 \sum_{i=0}^{m_j-1} h_i \tilde{u}_{ix\overline{t}}^2(k) = O(\tau),  \quad \text{as}~\tau \to 0.
\end{equation}

Let $\nu^{\tau}(t)=\nu^{k}, \ \mu^{\tau}(t)=\mu^{k}$ , if $t_{k-1}< t \leq t_{k}$, $k=1,\ldots,n$. We have
\begin{equation}
\Vert \nu^{\tau}-\nu \Vert_{L^{2}[0,T]}\to 0, \ \Vert \mu^{\tau}-\mu \Vert_{L^{2}[0,T]}\to 0~\text{as}~\tau \to 0\label{Eq:W:3:60}
\end{equation}
We estimate the first term in $\mathcal{I}_n(\mathcal{Q}_n(v))$ as follows
\begin{equation}
\beta_{0}\tau \sum_{j=1}^{n}|u_0(k)-\nu^{k}|^{2}=\beta_{0}\sum_{k=1}^{n}\int_{t_{k-1}}^{t_{k}}|u_0(k)-\nu^{k}|^{2}\, dt=\beta_{0}\int_{0}^{T}|u^{\tau}(0,t)-\nu^{\tau}(t)|^{2}\, dt\label{Eq:W:3:61}
\end{equation}
From (\ref{Eq:W:3:60}) it follows that
\begin{equation}
\lim_{n\to\infty}\beta_{0}\tau\sum_{k=1}^{n}|u_0(k)-\nu^{k}|^{2}=\beta_{0}\Vert u(0,t)-\nu(t) \Vert_{L^{2}[0,T]}^{2}\label{Eq:W:3:62}
\end{equation}
We estimate the second term in $\mathcal{I}_n(\mathcal{Q}_n(v))$ as follows
\begin{gather}
\beta_{1}\tau\sum_{k=1}^{n}|u_{m_k}(k)-\mu^{k}|^{2}=2\beta_{1}\sum_{k=1}^{n}\int_{t_{k-1}}^{t_{k}}\int_{s(t)}^{s_{k}}\frac{\partial u^{\tau}}{\partial x}\left(u^{\tau}(s(t),t)-\mu^{\tau}(t) \right)\, dx\, dt\nonumber\\ +\beta_{1}\sum_{k=1}^{n}\int_{t_{k-1}}^{t_{k}}|u(s(t);k)-\mu^{k}|^{2}\, dt  +\beta_{1} \sum_{k=1}^{n}\int_{t_{k-1}}^{t_{k}}\left( \int_{s(t)}^{s_{k}}\frac{\partial u^{\tau}}{\partial x}\, dx \right)^{2}\, dt=I_1+I_2+I_3\label{Eq:W:3:63}
\end{gather}
Since $\Big \Vert \frac{\partial u^\tau}{\partial x}\Big \Vert_{L_2(D)}$ and $\Vert u^\tau(s(t),t)-\mu^\tau \Vert_{L_2[0,T]}$ are uniformly bounded, and $\{s^n\}$ converges to $s$ uniformly on $[0,T]$, by applying CBS inequality and \eqref{alma} it easily follows that
\begin{equation}
\lim\limits_{n\to \infty} I_1 = 0, \ \lim\limits_{n\to \infty} I_3 = 0\label{Eq:W:3:65}
\end{equation}
It remains to prove that
\begin{equation}
\lim\limits_{n\to \infty} I_2=\lim\limits_{n\to \infty}\beta_{1}\int_{0}^{T}|u^{\tau}(s(t),t)-\mu^{\tau}(t)|^{2}\, dt = \beta_{1}\int_{0}^{T}|u(s(t),t)-\mu(t)|^{2}\, dt\label{Eq:W:3:64}
\end{equation}
 Since  $\mu^\tau \rightarrow \mu$ strongly in in $L_2[0,T]$ it is enough to show that
  \begin{gather}
    \Vert u^\tau(s(t),t) - u(s(t),t) \Vert_{L_2[0,T]} \ \to 0
    \intertext{as $\tau \to 0$. For any fixed $m>0$, estimate}
    \Vert u^\tau(s(t),t) - u(s(t),t) \Vert_{L_2[0,T]} \leq \Vert u^\tau(s(t),t) - u^\tau(s(t)-\epsilon_m,t)\Vert_{L_2[0,T]} + \nonumber
    \\
    +\Vert u^\tau(s(t)-\epsilon_m,t) - \hat{u}^\tau(s(t)-\epsilon_m,t)\Vert_{L_2[0,T]} + \nonumber
    \\
    + \Vert \hat{u}^\tau(s(t)-\epsilon_m,t) - u(s(t)-\epsilon_m,t)\Vert_{L_2[0,T]}
    + \Vert u(s(t)-\epsilon_m,t) - u(s(t),t)\Vert_{L_2[0,T]}\label{eq:InQnv-term2-est6-1}
    \intertext{Estimate the first term on the right-hand side of~\eqref{eq:InQnv-term2-est6-1} as}
    \Vert u^\tau(s(t),t) - u^\tau(s(t)-\epsilon_m,t)\Vert_{L_2[0,T]} = \left(\int_0^T \Big | \int_{s(t)-\epsilon_m}^{s(t)} \frac{\partial u^\tau(x,t)}{\partial x}\,dx \Big |^2 \,dt\right)^{1/2}
    \intertext{By CBS inequality and the first energy estimate}
    \Vert u^\tau(s(t),t) - u^\tau(s(t)-\epsilon_m,t)\Vert_{L_2[0,T]}
    \leq \sqrt{\epsilon_m}\Big \Vert \frac{\partial u^\tau}{\partial x}\Big \Vert_{L_2(D)} \leq C \sqrt{\epsilon_m}\label{eq:InQnv-term2-est6-2}
    \intertext{for $C$ independent of $n$.
      Similarly, the last term in~\eqref{eq:InQnv-term2-est6-1} is estimated by using CBS and energy estimate \eqref{V210estimate}:}
    \Vert u(s(t)-\epsilon_m,t) - u(s(t),t)\Vert_{L_2[0,T]}
    \leq C \sqrt{\epsilon_m}\label{eq:InQnv-term2-est6-3}
  \end{gather}
  Fix $\epsilon > 0$ and find $M$ such that for all $m \geq M$, $C \sqrt{\epsilon_m} \leq \epsilon/4$.
  Taking $m = M$, it follows from~\eqref{eq:InQnv-term2-est6-1}--\eqref{eq:InQnv-term2-est6-3}
  \begin{gather}
    \Vert u^\tau(s(t),t) - u(s(t),t)\Vert_{L_2[0,T]} \leq \frac{\epsilon}{2} 
    +\Vert u^\tau(s(t)-\epsilon_M,t) - \hat{u}^\tau(s(t)-\epsilon_M,t)\Vert_{L_2[0,T]}\nonumber\\
    + \Vert \hat{u}^\tau(s(t)-\epsilon_M,t) - u(s(t)-\epsilon_M,t)\Vert_{L_2[0,T]}\label{eq:InQnv-term2-est6-4}
  \end{gather}
  The second term in~\eqref{eq:InQnv-term2-est6-4} is estimated through Sobolev embedding of traces as
  \begin{gather}
    \Vert u^\tau(s(t)-\epsilon_M,t) - \hat{u}^\tau(s(t)-\epsilon_M,t)\Vert_{L_2[0,T]} \leq C \Vert \hat{u}^\tau - u^\tau \Vert_{W_2^{1,0}(\Omega_M)}
    \intertext{By \eqref{armud}, \eqref{Eq:W:3:59} there exists $\tau_0(M)>0$ such that $\forall \tau < \tau_0$}
   \Vert u^\tau(s(t)-\epsilon_M,t) - \hat{u}^\tau(s(t)-\epsilon_M,t)\Vert_{L_2[0,T]}  \leq \frac{\epsilon}{4}\label{garpiz}
  \end{gather}
  According to well-known compact embedding theorem weak convergence of $\hat{u}^\tau \to u$ in $W_2^{1,1}(\Omega_M)$ implies strong convergence of traces $\hat{u}^\tau|_{x=s(t)-\epsilon_M}$ to $u|_{x=s(t)-\epsilon_M}$ in $L_2[0,T]$ ~\cite{LSU}; that is, there exists $\tau_1(M)$ such that for all $\tau < \tau_1$ we have
  \begin{equation}\label{rrrr}
    \Vert \hat{u}^\tau(s(t)-\epsilon_M,t)- u(s(t)-\epsilon_M,t)\Vert_{L_2[0,T]} < \frac{\epsilon}{4}
  \end{equation}
Hence, according to \eqref{eq:InQnv-term2-est6-4}, \eqref{garpiz}, \eqref{rrrr}, for arbitrary $\epsilon>0$ we can find $\tau_2=\min(\tau_0; \tau_1)$ such that for all $\tau < \tau_2$ we have 
  \begin{equation}
  \Vert u^\tau(s(t),t)-u(s(t),t) \Vert_{L_2[0,T]} < \epsilon,
  \end{equation}
which proves \eqref{Eq:W:3:64}. Lemma is proved.

\begin{lemma}\label{condition(2)}
For arbitrary $[v]_{n}\in V_{R}^{n}$
\begin{equation}
\lim\limits_{n\to \infty} \Big ( \mathcal{J}(\mathcal{P}_n([v]_n))-\mathcal{I}_n([v]_n)  \Big )= 0\label{Eq:W:3:67}
\end{equation}
\end{lemma}
\textbf{Proof:} Let $[v]_{n}\in V_{R}^{n}$ and $v^{n}=(s^{n},g^{n})=\mathcal{P}_n([v]_n)$. From  Lemma~\ref{mappings}
it follows that the sequence $\{\mathcal{P}_n([v]_n\}$ is weakly precompact in $W_2^2[0,T]\times W_2^1[0,T]$. Assume that the whole sequence converges to $\tilde{v}=(\tilde{s},\tilde{g})$ weakly in $W_2^2[0,T]\times W_2^1[0,T]$. This implies the strong convegence in $W_2^1[0,T] \times L_2[0,T]$. From the well-known property of weak convergence it follows that $\tilde{v} \in V_{R}$. In particular $s^{n}$ converges to $\tilde{s}$ uniformly on $[0,T]$ and we have
\begin{equation}
\lim_{n\to\infty}\max_{0 \leq i \leq n}\left| s^{n}(t_{i})-\tilde{s}(t_{i})\right|=0\label{Eq:W:3:68}
\end{equation}
We have
\begin{equation}\label{Eq:W:3:69}
\mathcal{I}_{n}\big( [v]_{n}\big)-\mathcal{J}(v^{n}) = \mathcal{I}_{n}\big( [v]_{n}\big)-\mathcal{J}(\tilde{v}) + \mathcal{J}(\tilde{v})-\mathcal{J}(v^{n})
\end{equation}
Since $\mathcal{J}(v)$ is weakly continuous in $W_2^2[0,T] \times W_2^1[0,T]$ it follows that
\[
\lim_{n\to\infty}\left( \mathcal{J}(\tilde{v})-\mathcal{J}(v^{n})\right)=0. \]
Hence, we only need to prove that
\begin{equation}\label{Eq:W:3:70}
\lim_{n\to\infty} \mathcal{I}_{n}\big( [v]_{n}\big) = \mathcal{J}(\tilde{v})
\end{equation}
The proof of \eqref{Eq:W:3:70} is almost identical to the proof of Lemma~\ref{condition(1)}. Lemma is proved. 

Having Lemmas~\ref{condition(3)}, ~\ref{condition(1)} and ~\ref{condition(2)}, Theorem~\ref{convergence} 
follows from Lemma~\ref{generalcriteria}.



\end{document}

%% file: dibemod.0605


\newtheorem{theorem}{Theorem}[section]
\newtheorem{proposition}{Proposition}[section]
\newtheorem{lemma}{Lemma}[section]
\newtheorem{corollary}{Corollary}[section]
\newtheorem{remark}{Remark}[section]
\newtheorem{proof}{Proof:}

\renewcommand{\thesection}{\arabic{section}}
\renewcommand{\theequation}{\thesection.\arabic{equation}}
\renewcommand{\thetheorem}{\thesection.\arabic{theorem}}
\numberwithin{equation}{section}
\numberwithin{theorem}{section}
\numberwithin{proposition}{section}
\numberwithin{lemma}{section}
\numberwithin{remark}{section}
\setcounter{secnumdepth}{5}


\newcommand{\cl}{\centerline}
\newcommand{\sms}{\smallskip}
\newcommand{\ms}{\medskip}
\newcommand{\bs}{\bigskip}
\newcommand{\noi}{\noindent}
\newcommand{\itl}[1]{\textit{#1}}
\newcommand{\blf}[1]{\textbf{#1}}
\newcommand{\dsty}{\displaystyle}
\newcommand{\txty}{\textstyle}
\newcommand{\ssty}{\scriptstyle}
\newcommand{\tty}{\texttt}


\newcommand\Par{\mathhexbox278\,}


\newcommand{\al}{\alpha}
\newcommand{\Al}{\Alpha}
\newcommand{\be}{\beta}
\newcommand{\Be}{\Beta}
\newcommand{\Gm}{\Gamma}
\newcommand{\gm}{\gamma}
\newcommand{\dl}{\delta}
\newcommand{\Dl}{\Delta}
\newcommand{\lm}{\lambda}
\newcommand{\Lm}{\Lambda}
\newcommand{\kp}{\kappa}
\newcommand{\varep}{\varepsilon}
\newcommand{\vp}{\varphi}
\newcommand{\sig}{\sigma}
\newcommand{\Sig}{\Sigma}
\newcommand{\om}{\omega}
\newcommand{\Om}{\Omega}
\newcommand{\uom}{\mbox{\boldmath$\omega$}}
\newcommand{\btau}{\mbox{\boldmath$\tau$}}
\newcommand{\bnu}{\mbox{\boldmath$\nu$}}
\newcommand{\up}{\upsilon}
\newcommand{\z}{\zeta}


\newcommand{\df}[1]{\buildrel\mbox{\small def}\over{#1}}
\newcommand{\op}[1]{\buildrel\mbox{\tiny o}\over{#1}}
\newcommand{\db}{\prime\prime}
\newcommand{\bsl}{\backslash}
\newcommand{\lb}{\lbrack\!\lbrack}
\newcommand{\rb}{\rbrack\!\rbrack}
\newcommand\la{\langle}
\newcommand\ra{\rangle}
\newcommand{\ev}{\equiv}
\newcommand{\nev}{\not\equiv}
\newcommand{\nn}{\mathbb{N}}
\newcommand{\qq}{\mathbb{Q}}
\newcommand{\zz}{\mathbb{Z}}
\newcommand{\rr}{\mathbb{R}}
\newcommand{\rn}{\rr^N}
\newcommand{\cc}{\mathbb{C}}
\newcommand{\id}{\mathbb{I}}
\newcommand{\bo}{\mathbb{O}}

\newcommand{\amsb}[1]{\mathbb{#1}}
\newcommand{\mcl}[1]{\mathcal{#1}}
\newcommand{\bl}[1]{\mathbf{#1}}
\newcommand{\ov}[1]{\overline{#1}}
\newcommand{\wt}[1]{\widetilde{#1}}
\newcommand{\wh}[1]{\widehat{#1}}

\newcommand{\lra}{\longrightarrow}
\newcommand{\LLR}{\Longleftrightarrow}
\newcommand{\LRA}{\Longrightarrow}
\newcommand{\LLA}{\Longleftarrow}


\newcommand{\bbox}{\vrule height.6em width.6em 
depth0em} 
\newcommand{\os}{\vbox{\hrule \hbox{\vrule 
height.6em depth0pt 
\hskip.6em \vrule height.6em depth0em}
\hrule}} 


\newcommand{\dvg}{\operatorname{div}}
\newcommand{\curl}{\operatorname{curl}}
\newcommand{\supp}{\operatorname{supp}}
\newcommand{\essup}{\operatornamewithlimits{ess\,sup}}
\newcommand{\essinf}{\operatornamewithlimits{ess\,inf}}
\newcommand{\essosc}{\operatornamewithlimits{ess\,osc}}
\newcommand{\osc}{\operatornamewithlimits{osc}}
\newcommand{\sign}{\operatorname{sign}}
\newcommand{\loc}{\operatorname{loc}}
\newcommand{\diam}{\operatorname{diam}}
\newcommand{\dist}{\operatorname{dist}}
\newcommand{\card}{\operatorname{card}}
\newcommand{\meas}{\operatorname{meas}}
\newcommand{\spn}{\operatorname{span}}
\newcommand{\dtm}{\operatorname{det}}
%


\newcommand{\overlim}{\mathop{\overline{\lim}}\limits}
\newcommand{\underlim}{\mathop{\underline{\lim}}\limits}
\newcommand{\ttop}[2]{\genfrac{}{}{0pt}{}{#1}{#2}}
\newcommand{\bcu}{\mathop{\txty{\bigcup}}\limits}
\newcommand{\bca}{\mathop{\txty{\bigcap}}\limits}
\newcommand{\bsu}{\mathop{\txty{\sum}}\limits}
\newcommand{\pro}{\mathop{\txty{\prod}}\limits}


\newcommand{\pl}{\partial}
\newcommand{\ptt}{\frac{\pl}{\pl t}}
\newcommand{\ppx}{\frac\pl{\pl x}}
\newcommand{\dds}{\frac d{ds}}
\newcommand{\ddt}{\frac d{dt}}


\newcommand{\intl}{\int\limits}
\newcommand{\iintl}{\iint\limits}


\def\Xint#1{\mathchoice
    {\XXint\displaystyle\textstyle{#1}}%
    {\XXint\textstyle\scriptstyle{#1}}%
    {\XXint\scriptstyle\scriptscriptstyle{#1}}%
    {\XXint\scriptscriptstyle\scriptscriptstyle{#1}}%
    \!\int}
\def\XXint#1#2#3{\setbox0=\hbox{$#1{#2#3}{\int}$}
    \vcenter{\hbox{$#2#3$}}\kern-0.5\wd0}
\def\bint{\Xint-}
\def\dashint{\Xint{\raise4pt\hbox to7pt{\hrulefill}}}


\newcommand{\ovl}[3]{\int_{#1}^{#2}\kern-#3pt\raise4pt\hbox to7pt{\hrulefill}\ }

\newcommand{\ovll}[3]{\intl_{#1}^{#2}\kern-#3pt\raise4pt\hbox to7pt{\hrulefill}\ }

\newcommand{\tvl}[2]{\iint_{#1}\kern-#2pt\raise4pt\hbox to7pt{\hrulefill}\ }



\newcommand{\omt}{\Om_T}
\newcommand{\plo}{\partial\Omega}
\newcommand{\ovo}{\bar{\Om} }

%
\newcommand{\ci}[1]{C^\infty\!\left({#1}\right)}
\newcommand{\cio}[1]{C_o^\infty\!\left({#1}\right)}
\newcommand{\lloc}[1]{L_{\loc}\!\left({#1}\right)}
\newcommand{\xy}{|x-y|}


\newcommand{\intom}{\intl_{\Om}}
\newcommand{\intbo}{\intl_{\plo}}
\newcommand{\inom}{\int_{\Om}}
\newcommand{\inbo}{\int_{\plo}}
\newcommand{\intrn}{\intl_{\rn}}


\newcommand{\bye}{\end{document}}



%
%
%

%
%
%
%
%

%
%